# Stability Analysis of Nonlinear Inviscid Microscopic and Macroscopic Traffic Flow Models of Bidirectional Cruise-Controlled Vehicles

**Iasson Karafyllis[*], Dionysis Theodosis[**] and Markos Papageorgiou[**]**

[*]Dept. of Mathematics, National Technical University of Athens,
Zografou Campus, 15780, Athens, Greece,
emails: iasonkar@central.ntua.gr , iasonkaraf@gmail.com

[**] Dynamic Systems and Simulation Laboratory,
Technical University of Crete, Chania, 73100, Greece
(emails: dtheodosis@dssl.tuc.gr , markos@dssl.tuc.gr)

**Abstract**

The paper introduces a new bidirectional microscopic inviscid Adaptive Cruise Control (ACC) model that uses only spacing information from the preceding and following vehicles in order to select the proper control action to avoid collisions and maintain a desired speed. *KL* estimates that guarantee uniform convergence of the ACC model to the set of equilibria are provided. Moreover, the corresponding macroscopic model is derived, consisting of a conservation equation and a momentum equation that contains a nonlinear relaxation term. It is shown that, if the density is sufficiently small, then the macroscopic model has a solution that approaches exponentially the equilibrium speed (in the sup norm) while the density converges exponentially to a traveling wave. Numerical simulations are also provided, illustrating the properties of the microscopic and macroscopic inviscid ACC models.

## 1. Introduction

Microscopic traffic models describe the longitudinal (car-following) and lateral (lane-changing) movement of each single vehicle in the traffic stream, see [39]. Microscopic models based on Adaptive Cruise Control (ACC) and Cooperative Adaptive Cruise Control (CACC) systems are widely regarded as the basis of future generations of automated vehicles since they have the potential of increasing safety, reduce traffic accidents, and improve traffic flow on highways [18], [27]. Both ACC and CACC systems have been extensively studied in the literature (see for instance [6], [18], [22], [27], [33], [47] and references therein). The simplest form of interaction between vehicles, that gives more flexibility, than the typical Follow-the-Leader scheme [33], [39], is the bidirectional scheme which monitors the behavior of both the preceding and the following vehicles, see for instance [3], [11], [17], [25], [32], [41], [46]. It has been shown that bidirectional ACC systems can improve the platoon cohesion, [41], [47], and present strong attenuation of disturbances along a string of vehicles [4], [11], [14], [46].

Contrary to the microscopic models, macroscopic traffic models describe the traffic flow as a liquid that is characterized by macroscopic quantities, such as flow, density, and mean speed of vehicles. The first, simplest and most influential traffic flow model has been independently proposed by Lighthill and Whitham [26] and Richards [34] (LWR model). The well-known LWR model is a first-order model governed by one scalar hyperbolic partial differential equation that expresses the conservation of vehicles. First-order traffic models, while well-studied, present several limitations, as they do not capture significant traffic flow dynamics of interest [29]. To overcome the limitations of first-order models, second-order traffic flow models were considered, that introduce an additional differential equation for the mean speed, see for instance [1], [5], [10], [15], [28], [31], [43], [44], [40], [45]. Both first-order and second-order models have been extended to adjust the vehicle speed based on a perception of downstream density, see for instance [8], [9], [12], [37] and references therein. In the era of connected and automated vehicles, it is possible for vehicles to use backward sensors or to communicate their presence to other vehicles; hence, vehicles may adjust their speed based also on upstream density, in addition to downstream density. In [20], [30], the effect of the upstream density on the speed adjustment was termed as 'nudging', and it was shown that nudging can increase the flow in a ring-road and, if properly designed, can have a strong stabilizing effect on traffic flow.

The stability analysis of microscopic traffic flow models of vehicles under the effect of cruise controllers is challenging, because they are nonlinear models that evolve on sets which are not linear spaces. More specifically, the state space of microscopic traffic flow models of vehicles under ACC is usually neither an open set (see [38] and the recent paper [36] for the extension of the Input-to-State Stability property to systems evolving on open sets), nor a closed set. In addition, the set of equilibrium points can be an unbounded set and does not present uniform attractivity properties. In fact, the latter is a case for which very limited results exist in the literature (see [2], [7], [13], [34]). In this paper, we study a microscopic traffic flow model of vehicles under the effect of ACC which presents all the above characteristics.

This paper presents a novel bidirectional, microscopic, inviscid ACC model and its corresponding second-order macroscopic model. We call the model "inviscid" because it gives rise to a macroscopic model that contains no viscosity term. The proposed bidirectional microscopic model is based on the two-dimensional cruise controller for autonomous vehicles, recently proposed in [21], and we prove in this work that it has the following main features:
- Each vehicle uses only the distance from its preceding and following vehicles to select the proper control action (vehicle acceleration);
- the vehicles do not collide with each other;
- the speeds of all vehicles are always non-negative and remain below an a priori given speed limit;
- the ultimate distance between two consecutive vehicles is guaranteed to be greater than a pre-specified constant;
- all vehicle speeds converge to a given longitudinal speed set-point; and
- all the above features are valid globally, i.e., hold for all physically relevant initial conditions.

In addition to the above features, it is shown, by exploiting LaSalle's Invariance Principle (see [23]) that, the solutions of the microscopic inviscid model converge asymptotically to a set of equilibrium points from any (arbitrary) physically relevant initial condition (Theorem 2.1). However, since LaSalle's Invariance Principle does not guarantee a uniform attractivity property to the set of equilibria, we construct a strict Lyapunov function for the closed-loop system (Theorem 2.2). Using the constructed Lyapunov function, we establish a *KL* estimate for the solutions of the microscopic model that guarantees uniformity of the convergence rate to the set of equilibrium points (Theorem 2.3). We also show that, for specific initial conditions, the convergence is exponential (Proposition 3.1). The main theoretical challenges stem from the fact that the control system studied in the paper evolves on a specific set (which is neither open nor closed), and, in addition, various objectives and constraints must be satisfied simultaneously and globally (positive speeds within road speed limits that converge to a specific speed set-point). In [21], it was not possible to show that vehicles operating on the two-dimensional surface of a lane-free road ultimately reach a set of equilibrium configurations. However, for the case of a string of vehicles moving only longitudinally, we prove that the vehicles do reach a set of equilibrium configurations, where the distance between two consecutive vehicles is guaranteed to be bounded and be greater than a pre-specified constant. Furthermore, we examine by simulation the speed and spacing first-to-last amplification factor of the inviscid microscopic model (Section 5). Such metrics have been widely used to study the robustness of bidirectional ACC models to external disturbances (see for instance [4], [14], [17]). We observe that disturbances dissipate much faster along the string of vehicles for the bidirectional inviscid model than the collision-free Follow-the-Leader model proposed in [22]. Therefore, our simulation experiments show that nonlinear controllers can exhibit string stability properties, contrary to linear approaches that require non-identical controllers (see [24]). Finally, we (formally) derive the macroscopic inviscid model that corresponds to the bidirectional microscopic model, and show that, if the initial density is sufficiently small, then the macroscopic model has a solution that approaches exponentially the equilibrium speed (in the sup norm) while the density converges exponentially to a traveling wave (Theorem 3.2).

The structure of the paper is as follows. Section 2 is devoted to the presentation of the bidirectional microscopic inviscid ACC model and its stability properties. Section 3 presents the corresponding macroscopic inviscid model and its analogy to the microscopic model. Section 4 presents numerical examples to demonstrate the properties of the microscopic and macroscopic inviscid models. In Section 5, we investigate by simulation the first-to-last amplification factor and the sensitivity to external disturbances of the microscopic model. All proofs of the main results are provided in Section 6. Finally, some concluding remarks are given in Section 7.

**Notation.** Throughout this paper, we adopt the following notation.

* $\Re_+ := [0, +\infty)$ denotes the set of non-negative real numbers.

* By $|x|$ we denote both the Euclidean norm of a vector $x \in \Re^n$ and the absolute value of a scalar $x \in \Re$.

* By $K$ we denote the class of increasing $C^0$ functions $a: \Re_+ \to \Re_+$ with $a(0) = 0$. By $K_\infty$ we denote the class of increasing $C^0$ functions $a: \Re_+ \to \Re_+$ with $a(0) = 0$ and $\lim_{s \to +\infty} a(s) = +\infty$. By $KL$ we denote the set of all continuous functions $\sigma: \Re_+ \times \Re_+ \to \Re_+$ with the properties: (i) for each $t \geq 0$ the mapping $\sigma(\cdot, t)$ is of class $K$; (ii) for each $s \geq 0$, the mapping $\sigma(s, \cdot)$ is non-increasing with $\lim_{t \to +\infty} \sigma(s, t) = 0$.

* By $C^0(A, \Omega)$, we denote the class of continuous functions on $A \subseteq \Re^n$, which take values in $\Omega \subseteq \Re^m$. By $C^k(A; \Omega)$, where $k \geq 1$ is an integer, we denote the class of functions on $A \subseteq \Re^n$ with continuous derivatives of order $k$, which take values in $\Omega \subseteq \Re^m$. When $\Omega = \Re$ the we write $C^0(A)$ or $C^k(A)$.

* Let $I \subseteq \Re$ be a given interval. $L^\infty(I)$ denotes the set of equivalence classes of measurable functions $f: I \to \Re$ for which $\|f\|_\infty = \operatorname{ess\,sup}_{x \in I}(|f(x)|) < +\infty$. By $W^{k,\infty}(I)$, where $k \geq 1$ is an integer, we denote the Sobolev spaces of functions $f \in L^\infty(I)$ which have weak derivatives of order $\leq k$, all of which belong to $L^\infty(I)$.

* For a set $S \subseteq \Re^n$, $\overline{S}$ denotes the closure of $S$.

* We denote by $dist(x, A)$ the Euclidean distance of the point $x \in \Re^n$ from the set $A \subset \Re^n$, i.e., $dist(x, A) = \inf\{|x - y| : y \in A\}$.

* Let $u : \mathfrak{R}_+ \times \mathfrak{R} \to \mathfrak{R}$, $(t,x) \to u(t,x)$ be any function differentiable with respect to its arguments. We use the notation $u_t(t,x) = \frac{\partial u}{\partial t}(t,x)$ and $u_x(t,x) = \frac{\partial u}{\partial x}(t,x)$ for the partial derivatives of $u$ with respect to $t$ and $x$, respectively. We use the notation $u[t]$ to denote the profile at certain $t \geq 0$, $(u[t])[x] := u(t,x)$, for all $x \in \mathfrak{R}$.

## 2. The Microscopic Inviscid ACC Model

### *2.1. Description of the model*

The movement of $n$ identical vehicles on a straight road under the cruise controller that was proposed in [21], when the vehicles are constrained to move on a line (longitudinal motion), is described by the following set of ODEs:

$$\begin{aligned}
\dot{x}_i &= v_i \quad, \quad i = 1, 2, ..., n \\
\dot{v}_1 &= -k_1(s_2)(v_1 - v^*) - V'(s_2) \\
\dot{v}_i &= -k_i(s_i, s_{i+1})(v_i - v^*) + V'(s_i) - V'(s_{i+1}) \quad, \quad i = 2, ..., n-1 \\
\dot{v}_n &= -k_n(s_n)(v_n - v^*) + V'(s_n)
\end{aligned} \quad (2.1)$$

where

$$\begin{aligned}
s_i &= x_{i-1} - x_i \quad, \quad i = 2, ..., n \\
k_1(s_2) &= \mu + g(-V'(s_2)) \\
k_i(s_i, s_{i+1}) &= \mu + g(V'(s_i) - V'(s_{i+1})) \quad, \quad i = 2, ..., n-1 \\
k_n(s_n) &= \mu + g(V'(s_n))
\end{aligned} \quad (2.2)$$

$\mu, v^* > 0$ are constants, $V \in C^2((L, +\infty); \mathfrak{R}_+)$ is a potential function that satisfies

$$\begin{aligned}
&\lim_{s \to L^+}(V(s)) = +\infty, \\
&V''(s) \geq 0, \\
&s \geq \lambda \Rightarrow V(s) = 0, \\
&L < s < \lambda \Rightarrow V'(s) < 0
\end{aligned} \quad (2.3)$$

where $\lambda > L > 0$ are constants, and

$$g(s) = \frac{v_{max} f(s)}{v^*(v_{max} - v^*)} - \frac{s}{v^*}, \quad (2.4)$$

$v_{max} > v^*$ is a constant (the road speed limit), and $f \in C^1(\mathfrak{R})$ is a non-decreasing function that satisfies

$$\max(x, 0) \leq f(x), \text{ for all } x \in \mathfrak{R}. \quad (2.5)$$

Using (2.2), the model can be written in the following form

$$\dot{s}_i = v_{i-1} - v_i \quad , \quad i = 2,...,n$$
$$\dot{v}_1 = -k_1(s_2)\left(v_1 - v^*\right) - V'(s_2)$$
$$\dot{v}_i = -k_i(s_i, s_{i+1})\left(v_i - v^*\right) + V'(s_i) - V'(s_{i+1}) \quad , \quad i = 2,...,n-1 \tag{2.6}$$
$$\dot{v}_n = -k_n(s_n)\left(v_n - v^*\right) + V'(s_n)$$

Due to various constraints, such as minimum inter-vehicle distance and speeds within certain speed limits, the state space of model (2.6) is

$$\Omega = \left\{ (s_2,...,s_n, v_1,...,v_n) \in \Re^{2n-1} : \min_{i=2,...,n}(s_i) > L, \max_{i=1,...,n}(v_i) \leq v_{\max}, \min_{i=1,...,n}(v_i) \geq 0 \right\} \tag{2.7}$$

where $L$ is a given positive constant (the minimum distance between two vehicles for which the vehicles do not collide with each other). In what follows, we refer to model (2.6) as the "microscopic, inviscid ACC model", since the macroscopic analogue of (2.6) does not contain a viscosity term (see Section 3). Clearly, model (2.6) is nonlinear not only because of the nonlinearities appearing in the right-hand sides of (2.6), but also due to the fact that the state space $\Omega$ is not a linear subspace of $\Re^{2n-1}$. It should be also noticed that the state space is not an open set (see the recent paper [36] for the extension of the Input-to-State Stability property to systems defined on open sets) and it is not a closed set.

*2.2. Stability Analysis*

Since $v^* \in (0, v_{\max})$, it follows by (2.3) and (2.6), that the set

$$S = \left\{ (s_2,...,s_n, v_1,...,v_n) \in \Re^{2n-1} : \min_{i=2,...,n}(s_i) \geq \lambda, v_i = v^*, i = 1,...,n \right\} \subset \Omega \tag{2.8}$$

is the set of equilibrium points for the model. In what follows, we use the notation

$$s = (s_2,...,s_n) \in \Re^{n-1}, \quad v = (v_1,...,v_n) \in \Re^n. \tag{2.9}$$

Moreover, in what follows, we omit the arguments of the functions $k_i, i=1,...,n$, defined by (2.2) (for simplicity).

In this section we show the stability properties of the invariant set $S$. However, before we proceed, it is necessary to state clearly the characteristics of the problem that indicate the challenge (from a mathematical point of view) of the performed stability analysis:
 (i) system (2.6) is nonlinear,
 (ii) the state space of system (2.6) is neither a closed set, nor an open set, and
 (iii) the invariant set, whose stability properties are to be investigated, is not a bounded set (and consequently not a compact set).

We define the function $H: \Omega \to \Re_+$ by the formula

$$H(s,v) = \frac{1}{2}\sum_{i=1}^{n}\left(v_i - v^*\right)^2 + \sum_{i=2}^{n} V(s_i). \tag{2.10}$$

Notice that the function $H$ is nothing else but the mechanical energy of the system of $n$ vehicles relative to an observer that moves with constant speed $v^*$. Using (2.6) and (2.10), we obtain the following equation for all $(s,v) \in \Omega$ for the time derivative of the function $H$ along the solutions of (2.6):

$$\dot{H}(s,v) = -\sum_{i=1}^{n} k_i \left(v_i - v^*\right)^2. \tag{2.11}$$

It should be noticed that the function $H$ is not a strict Lyapunov function, because (2.11) shows that the derivative of $H$ can be zero for points out of the invariant set $S$. However, using (2.10) and (2.11), we can prove the following theorem.

**Theorem 2.1:** *For every initial condition $(s(0), v(0)) \in \Omega$, the solution $(s(t), v(t)) \in \Omega$ of (2.6) is defined for all $t \geq 0$ and satisfies*

$$s_i(t) \leq \max\left(\lambda, s_i(0)\right) + \mu^{-1} v_{\max}, \text{ for all } t \geq 0 \text{ and } i = 2,...,n. \tag{2.12}$$

*Moreover, $\lim_{t \to +\infty} (v_i(t)) = v^*$ for all $i = 1,...,n$ and $\lim_{t \to +\infty} (V(s_i(t))) = 0$ for all $i = 2,...,n$.*

It should be noted that Theorem 2.1 is not a specialization of Theorem 1 in [21] for vehicles moving on a straight line. Indeed, an application of Theorem 1 in [21] does not guarantee that the distance of two consecutive vehicles (i.e., $s_i$) is bounded from above. In other words, an application of Theorem 1 in [21] does not give us estimates (2.12). Notice also that, since $\lim_{t \to +\infty} (V(s_i(t))) = 0$ for all $i = 2,...,n$, properties (2.3) guarantee that $\liminf_{t \to +\infty} (s_i(t)) \geq \lambda$ for all $i = 2,...,n$, and that $\lim_{t \to +\infty} \left(dist((s(t), v(t)), S)\right) = 0$ (recall definition (2.8)). The proof of Theorem 2.1 relies on LaSalle's principle. However, LaSalle's principle does not guarantee *uniform* attraction to the set $S$. In order to be able to show uniform global attractivity properties for the set $S$, we need to provide a strict Lyapunov function for system (2.6). This is done by the following theorem.

**Theorem 2.2:** *For every $\beta > 0$, there exist non-decreasing functions $R \in C^1(\mathfrak{R}_+; (0, +\infty))$, $\kappa \in C^0(\mathfrak{R}_+; (0, +\infty))$ such that the following inequalities hold for all $(s,v) \in \Omega$:*

$$H(s,v) \leq W(s,v) \leq \kappa(H(s,v)) H(s,v) \tag{2.13}$$

$$\dot{W}(s,v) \leq -\beta \mu \sum_{i=1}^{n} \left(v_i - v^*\right)^2 - \frac{1}{8} \sum_{i=2}^{n} 4^i \left(V'(s_i)\right)^2 \tag{2.14}$$

*where $W : \Omega \to \mathfrak{R}_+$ is defined by the equation*

$$W(s,v) := R(H(s,v)) H(s,v) - \sum_{i=2}^{n} 4^i V'(s_i) \left(v_i - v^*\right), \text{ for all } (s,v) \in \Omega \tag{2.15}$$

*and $\dot{W}(s,v)$ denotes the time derivative of $W$ along the solutions of (2.6).*

**Remark:** Notice that $W(s,v) > 0$, $\dot{W}(s,v) < 0$ when $(s,v) \in \Omega \setminus S$, and $W(s,v) = \dot{W}(s,v) = 0$ when $(s,v) \in S$. Thus, the function $W$ defined by (2.15) is a strict Lyapunov function for the microscopic inviscid ACC model (2.6).

Using Theorem 2.2 we are in a position to prove that a $KL$ estimate holds for the solutions of (2.6). This estimate is important, because it guarantees uniformity of the convergence rate to the set $S$, as well as useful robustness properties (see the discussion in [38]; uniform global asymptotic stability with respect to two measures).

**Theorem 2.3:** *There exist a function $\sigma \in KL$ and a function $a \in K_\infty$ such that for every initial condition $(s(0), v(0)) \in \Omega$ the solution $(s(t), v(t)) \in \Omega$ of (2.6) is defined for all $t \geq 0$ and satisfies*

$$a\big(dist\big((s(t), v(t)), S\big)\big) \leq W(s(t), v(t)) \leq \sigma\big(W(s(0), v(0)), t\big), \text{ for all } t \geq 0. \tag{2.16}$$

**Remarks: (i)** While estimate (2.16) implies a uniform rate of convergence to the invariant set $S$, it does not imply an exponential rate of convergence. This fact can also be seen in the simulation examples of Section 5, see Example 1. A special case for exponential convergence to the set $S$ and its relation to the macroscopic inviscid model that corresponds to model (2.6) will be discussed in the following section (Proposition 3.1).
**(ii)** Theorem 2.1, Theorem 2.2 and Theorem 2.3 show global convergence of the solutions of the model (2.6) to the invariant set of equilibrium points $S$. However, this does not mean that every solution of (2.6) converges to an equilibrium point. In fact, we cannot conclude that the limits $\lim_{t \to +\infty} (s_i(t))$ exist for $i = 2, ..., n$. However, Theorem 2.1 shows that every solution of (2.6) satisfies the following estimates for $i = 2, ..., n$:

$$\lambda \leq \liminf_{t \to +\infty} (s_i(t)) \leq \limsup_{t \to +\infty} (s_i(t)) \leq \max(\lambda, s_i(0)) + \mu^{-1} v_{max}. \tag{2.17}$$

**(iii)** It should be noted that the stability estimate (2.16) does not establish Uniform Global Asymptotic Stability of the invariant set $S$, i.e., we do not show an estimate of the form $dist\big((s(t), v(t)), S\big) \leq \bar{\sigma}\big(dist\big((s(0), v(0)), S\big), t\big)$, for all $t \geq 0$ for a $KL$ function $\bar{\sigma}$. Instead, the stability estimate (2.16) establishes Uniform Global Asymptotic Stability with respect to the measures $\omega_1(s, v) = dist\big((s, v), S\big)$ and $\omega_2(s, v) = W(s, v)$ (see [38]).

## 3. The Macroscopic Inviscid ACC Model

In this section we focus on the macroscopic traffic model that corresponds to the microscopic model (2.1). Macroscopic traffic models describe the collective behavior of vehicles on a highway and allow for more direct insights and mathematical analysis of traffic properties, such as capacity flow and mean speed of a traffic stream, see for instance [1], [12], [26], [29], [42]. Various approaches have been suggested to derive macroscopic models for conventional traffic from microscopic models, see for instance [1], [9], [15], [16], [31], [45], and references therein.

### *3.1. The PDE model*

Let $\rho_{max}, v_{max} > 0$ and $v^* \in (0, v_{max})$, $\bar{\rho} \in (0, \rho_{max})$ be constants and let $\Phi : (0, \rho_{max}) \to \Re_+$ be a $C^3\big((0, \rho_{max})\big)$ non-negative function that satisfies:

$$\lim_{\rho \to \rho_{max}^-} (\Phi(\rho)) = +\infty, \quad \Phi(\rho) = 0 \text{ for all } \rho \in (0, \bar{\rho}] \tag{3.1}$$

$$\Phi'(\rho) > 0, \text{ for all } \rho \in (\bar{\rho}, \rho_{max}) \tag{3.2}$$

$$\Phi''(\rho) \geq 0, \text{ for all } \rho \in (0, \rho_{max}) \qquad (3.3)$$

The macroscopic model that corresponds to the microscopic model (2.1), as the number of vehicles $n$ tends to infinity and the potential function is given by $V(s) = \Phi\left(\dfrac{m}{ns}\right)$ with $\dfrac{m}{n}$ being the mass of every single vehicle, is formally derived in Section 6 and is the following nonlinear system of PDEs for $(t, x) \in (0, +\infty) \times \Re$

$$\rho_t(t,x) + v(t,x)\rho_x(t,x) + \rho(t,x)v_x(t,x) = 0$$
$$v_t(t,x) + v(t,x)v_x(t,x) - \Xi(t,x) = -\big(\mu + g(\Xi(t,x))\big)\big(v(t,x) - v^*\big) \qquad (3.4)$$

where $\rho(t, x)$ is the traffic density, $v(t, x)$ is the mean speed, and
$$\Xi(t,x) := -\rho_x(t,x)\big(2\Phi'(\rho(t,x)) + \rho(t,x)\Phi''(\rho(t,x))\big) \qquad (3.5)$$

with constrained values
$$0 < \rho(t,x) < \rho_{max} \; , \; 0 \leq v(t,x) \leq v_{max}, \text{ for all } (t,x) \in \Re_+ \times \Re. \qquad (3.6)$$

Four things should be noted about the nonlinear model (3.4)-(3.5):
1) there are no non-local terms in the model, despite the fact that the cruise controller proposed in [21] induces "nudging";
2) there are infinite equilibrium points for the model, namely the points where $v \equiv v^*$ and $\rho(x) \leq \bar{\rho}$ for all $x \in \Re$;
3) it is a second-order model; and
4) the model is highly nonlinear (it is not semilinear as the ARZ model [1], [45] or the PW model [31], [40]), due to the presence of a highly nonlinear relaxation term $-\big(\mu + g(\Xi(t,x))\big)\big(v(t,x) - v^*\big)$ in the speed PDE.

### *3.2. An analogy between the microscopic and the macroscopic model*

For the microscopic model (2.6), the following statement holds: When the vehicles have large initial distances between them, then they simply adjust their speeds without affecting each other. This is shown by the following proposition, whose proof is very simple and is omitted.

**Proposition 3.1:** *Suppose that* $s_i(0) \geq \max\big(\lambda - \omega^{-1}(v_{i-1}(0) - v_i(0)), \lambda\big)$ *for* $i = 2,...,n$, *where* $\omega = \mu + g(0)$. *Then the solution of the model (2.6) is given by the equations:*

$$v_i(t) = v^* + \exp(-\omega t)\big(v_i(0) - v^*\big) \; , \; i = 1,...,n$$
$$s_i(t) = s_i(0) + \omega^{-1}\big(v_{i-1}(0) - v_i(0)\big)\big(1 - \exp(-\omega t)\big) \; , \; i = 2,...,n \qquad (3.7)$$

In this case we have exponential convergence to the set $S$. This is also illustrated in the numerical examples of Section 5.
Similarly with Proposition 3.1, if the initial density is sufficiently small, then the macroscopic model (3.4)-(3.5) has a solution that approaches the equilibrium speed (in the sup norm), while the density remains small. The following theorem guarantees this fact.

**Theorem 3.2:** *Consider the initial-value problem*

$$\begin{aligned} \rho_t + v\rho_x + \rho v_x &= 0 \\ v_t + vv_x &= -\omega(v - v^*) \end{aligned} \quad \text{for } t \geq 0, \; x \in \Re \qquad (3.8)$$

$$\begin{aligned} \rho(0, x) &= \rho_0(x) \\ v(0, x) &= v_0(x) \end{aligned}, \text{ for } x \in \Re \qquad (3.9)$$

*where* $\omega = \mu + g(0)$, $\rho_0 \in C^1(\Re) \cap W^{1,\infty}(\Re), v_0 \in C^2(\Re) \cap W^{2,\infty}(\Re)$ *with* $\inf_{x \in \Re}(v_0'(x)) > -\omega$ *and* $\rho_0(x) > 0$ *for all* $x \in \Re$.

*Then, the initial-value problem (3.8), (3.9) has a unique solution that satisfies the estimates:*

$$\sup_{x \in \Re}(\rho(t, x)) \leq \frac{\omega \sup_{x \in \Re}(\rho_0(x))}{\omega + (1 - \exp(-\omega t)) \inf_{x \in \Re}(v_0'(x))}, \text{ for all } t \geq 0 \qquad (3.10)$$

$$\sup_{x \in \Re}(|v(t, x) - v^*|) \leq \exp(-\omega t) \sup_{x \in \Re}(|v_0(x) - v^*|), \text{ for all } t \geq 0 \qquad (3.11)$$

$$\rho(t, x) > 0, \text{ for all } t \geq 0, \; x \in \Re. \qquad (3.12)$$

*Moreover, there exists a function* $f : \Re \to (0, +\infty)$ *of class* $C^1(\Re) \cap L^\infty(\Re)$ *for which the following estimate holds:*

$$\sup_{x \in \Re, t \geq 0} \left(|\rho(t, x) - f(x - v^* t)| \exp(\omega t)\right) < +\infty. \qquad (3.13)$$

By virtue of (3.1) and (3.5), it follows that $\Xi(t, x) \equiv 0$ when $\sup_{x \in \Re}(\rho(t, x)) \leq \bar{\rho}$ for all $t \geq 0$. Therefore, in this case the PDE model (3.4)-(3.5) becomes identical to the PDE model (3.8). Consequently, Theorem 3.2 guarantees that, if the initial conditions satisfy the requirements

$$\inf_{x \in \Re}(v_0'(x)) > -\omega \qquad (3.14)$$

$$\sup_{x \in \Re}(\rho_0(x)) \leq \bar{\rho}\left(1 + \omega^{-1} \min\left(0, \inf_{x \in \Re}(v_0'(x))\right)\right) \qquad (3.15)$$

then the macroscopic model (3.4)-(3.5) has a solution that satisfies the estimates

$$\sup_{x \in \Re}(\rho(t, x)) \leq \bar{\rho}, \text{ for all } t \geq 0 \qquad (3.16)$$

$$\sup_{x \in \Re}(|v(t, x) - v^*|) \leq \exp(-\omega t) \sup_{x \in \Re}(|v_0(x) - v^*|), \text{ for all } t \geq 0 \qquad (3.17)$$

$$\rho(t, x) > 0, \text{ for all } t \geq 0, \; x \in \Re \qquad (3.18)$$

and there exists a function $f : \Re \to (0, \bar{\rho}]$ of class $C^1(\Re) \cap L^\infty(\Re)$ for which the following estimate holds:

$$\sup_{x \in \Re, t \geq 0} \left(|\rho(t, x) - f(x - v^* t)| \exp(\omega t)\right) < +\infty$$

Notice that the solution converges exponentially (in the sup norm) to the set of equilibrium points of the macroscopic model (but not necessarily to one equilibrium point). This is due to the fact that the whole profile ultimately moves with speed $v^*$ and thus we have a traveling wave.

## 4. Illustrative Examples

In this section we illustrate by simulations the properties of the microscopic and macroscopic inviscid models (2.6) and (3.8) respectively.

**Example 1:** (Asymptotic stability of the set $S$). In this example, we illustrate the results of Theorem 2.1, and demonstrate the convergence of the Lyapunov function $H(s,v)$ and the asymptotic convergence of the spacing and speed to the set of equilibrium points $S$ in (2.8). Consider the inviscid model (2.6) and let $V$ defined by

$$V(q) = \begin{cases} \dfrac{(\lambda-q)^3}{q-L}, & L < q \le \lambda \\ 0 & q > \lambda \end{cases} \quad (4.1)$$

Moreover, we define the function $f$ in (2.4) by means of

$$f(x) = \dfrac{1}{2\varepsilon} \begin{cases} 0 & \text{if } x \le -\varepsilon \\ (x+\varepsilon)^2 & \text{if } -\varepsilon < x < 0 \\ \varepsilon^2 + 2\varepsilon x & \text{if } x \ge 0 \end{cases} \quad (4.2)$$

which satisfies inequality (2.5) for any $\varepsilon > 0$. We consider $n=6$ vehicles with initial spacing and speed $s_i(0) \in (16,24)$, $v_i(0) \in (27,34)$, respectively. We also let $L = 5m$, $\lambda = 20m$, $v^* = 30m/s$, $v_{\max} = 35m/s$, $\varepsilon = 0.2$, and $\mu = 0.5$. The inter-vehicle distances and the speed of each vehicle are shown in the Figure 1. It is seen that the speeds of the vehicles converge to the set-point $v^*$, i.e., that $\lim_{t \to +\infty} (v_i(t)) = v^*$; and that the distances satisfy $\liminf_{t \to +\infty} (s_i(t)) \ge \lambda$. Finally, Figure 2 shows the evolution of the Lyapunov function $H(s,v)$, defined by (2.10), and its logarithm, indicating asymptotic convergence, but not exponential convergence (because the plot of $\ln(H(s,v))$ with respect to $t$ is not below a straight line of negative slope).

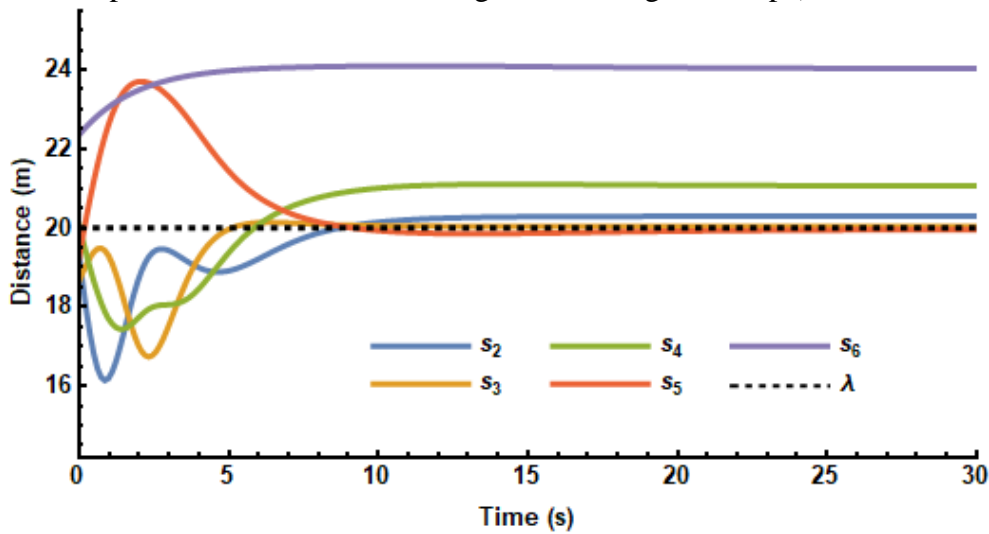

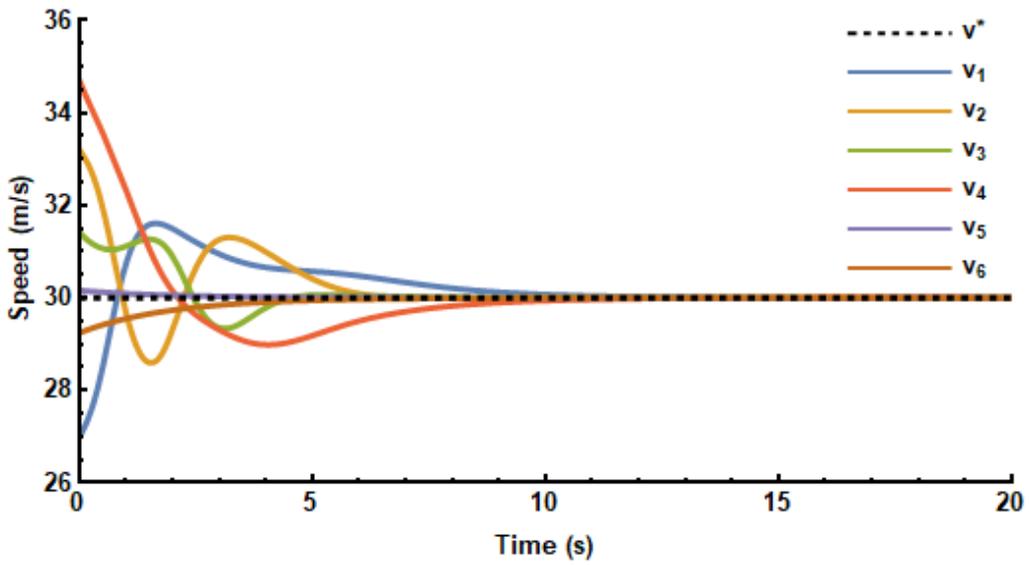

**Figure 1:** Inter-vehicle distances (top) and vehicle speeds (bottom).

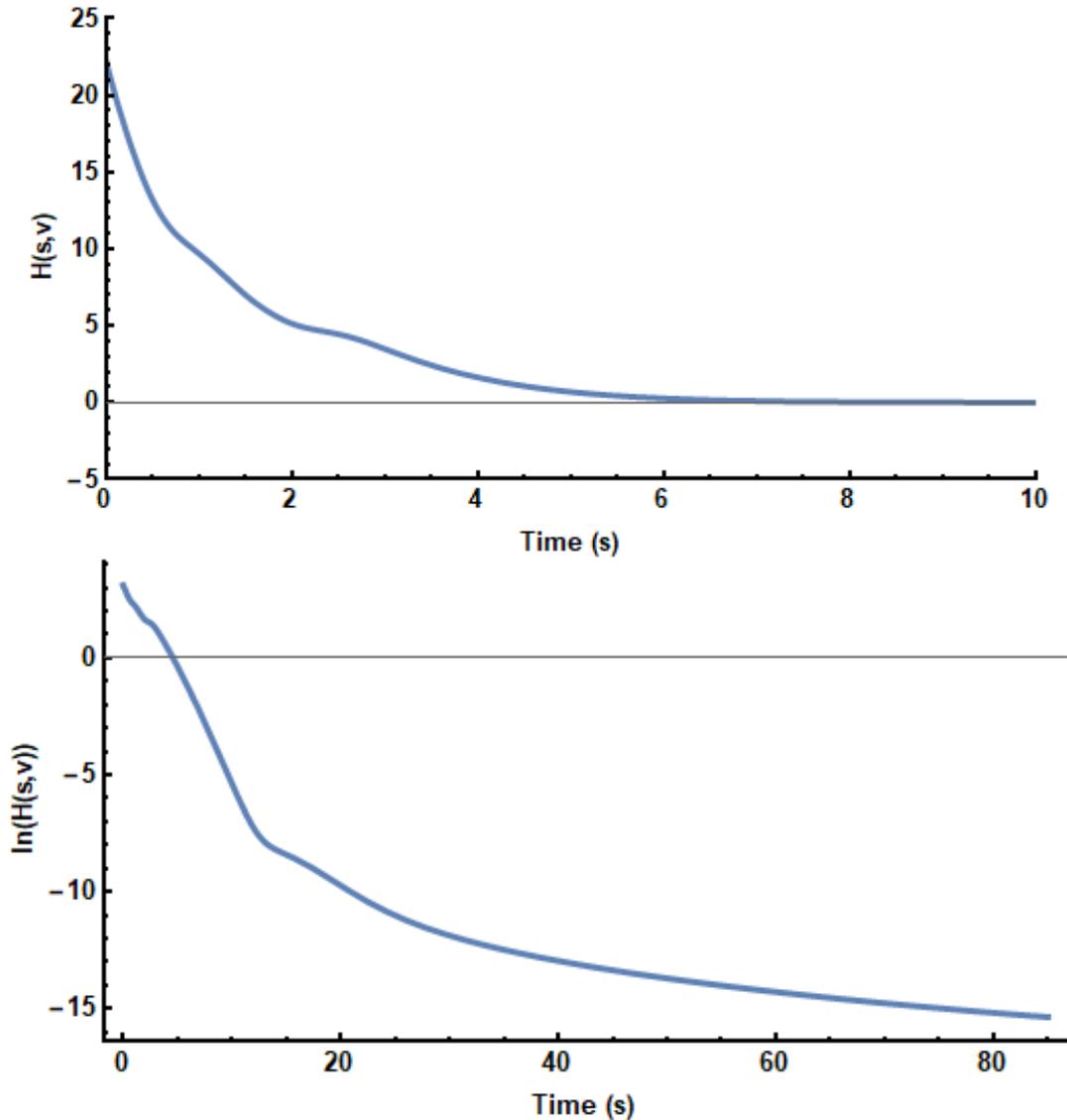

**Figure 2:** Convergence of the function $H(s,v)$ (top) and the graph of $\ln(H(s,v))$ (bottom) indicating asysmptotic convergence but not exponential convergence.

**Example 2:** (Exponential convergence to the set $S$) As indicated by Proposition 3.1, for the exponential convergence to the set of equilibrium points $S$, the initial spacing needs to satisfy the condition $s_i(0) \geq \max\left(\lambda - \omega^{-1}(v_{i-1}(0) - v_i(0)), \lambda\right)$ for each $i = 1,...,n$, where $\omega = \mu + g(0)$ and $g$ is given by (2.4). Values of $L = 5m$, $\lambda = 20m$, $v^* = 30m/s$, $v_{max} = 35m/s$, and $\mu = 0.5$ were used. For initial conditions satisfying this constraint, we have the speeds and inter-vehicle distances shown in Figure 3. The exponential convergence to the set $S$ is demonstrated by Figure 4, which shows the evolution of the Lyapunov function $H(s,v)$ and its logarithm $\ln(H(s,v))$.

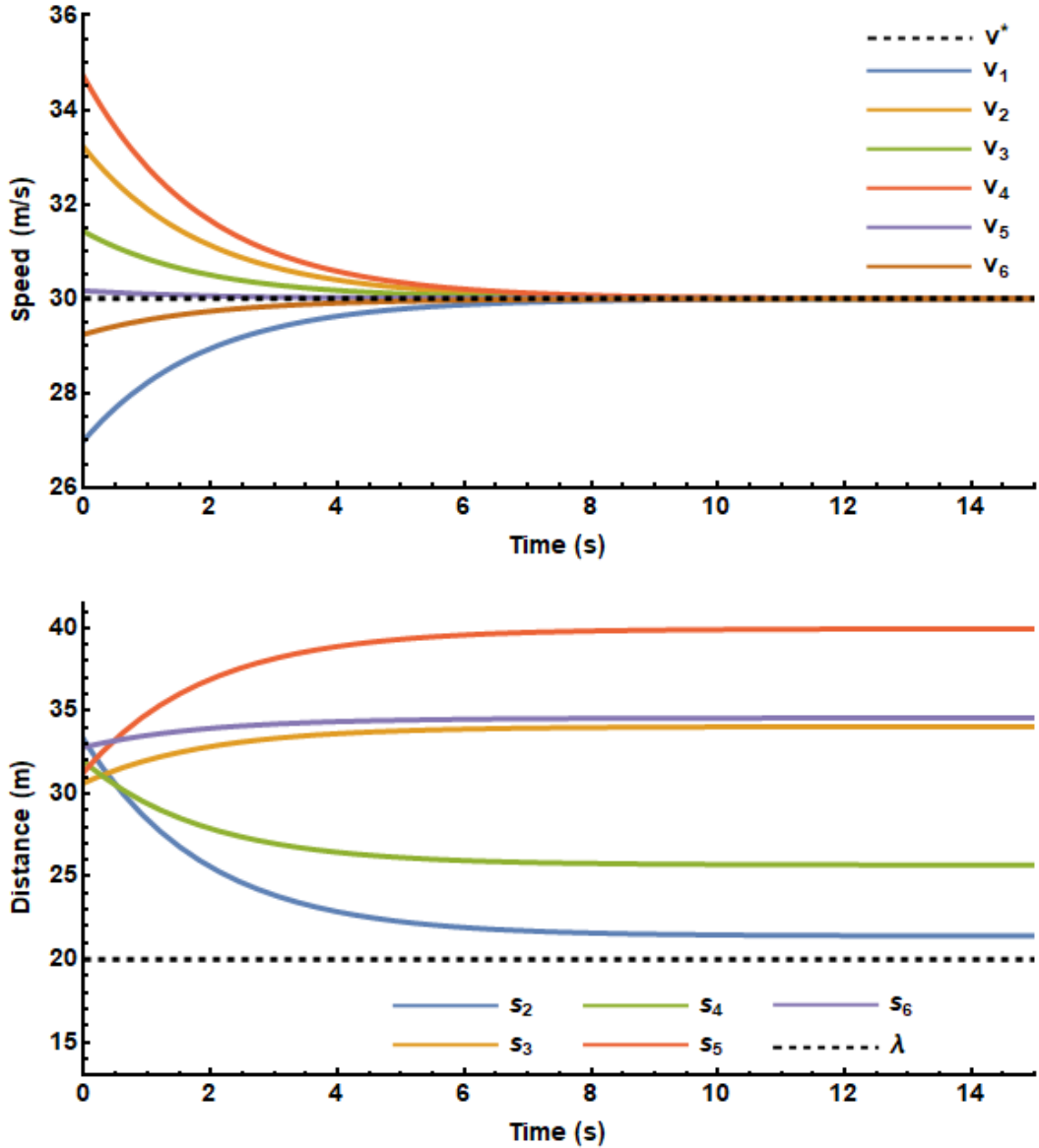

**Figure 3:** Exponential convergence to the set $S$.

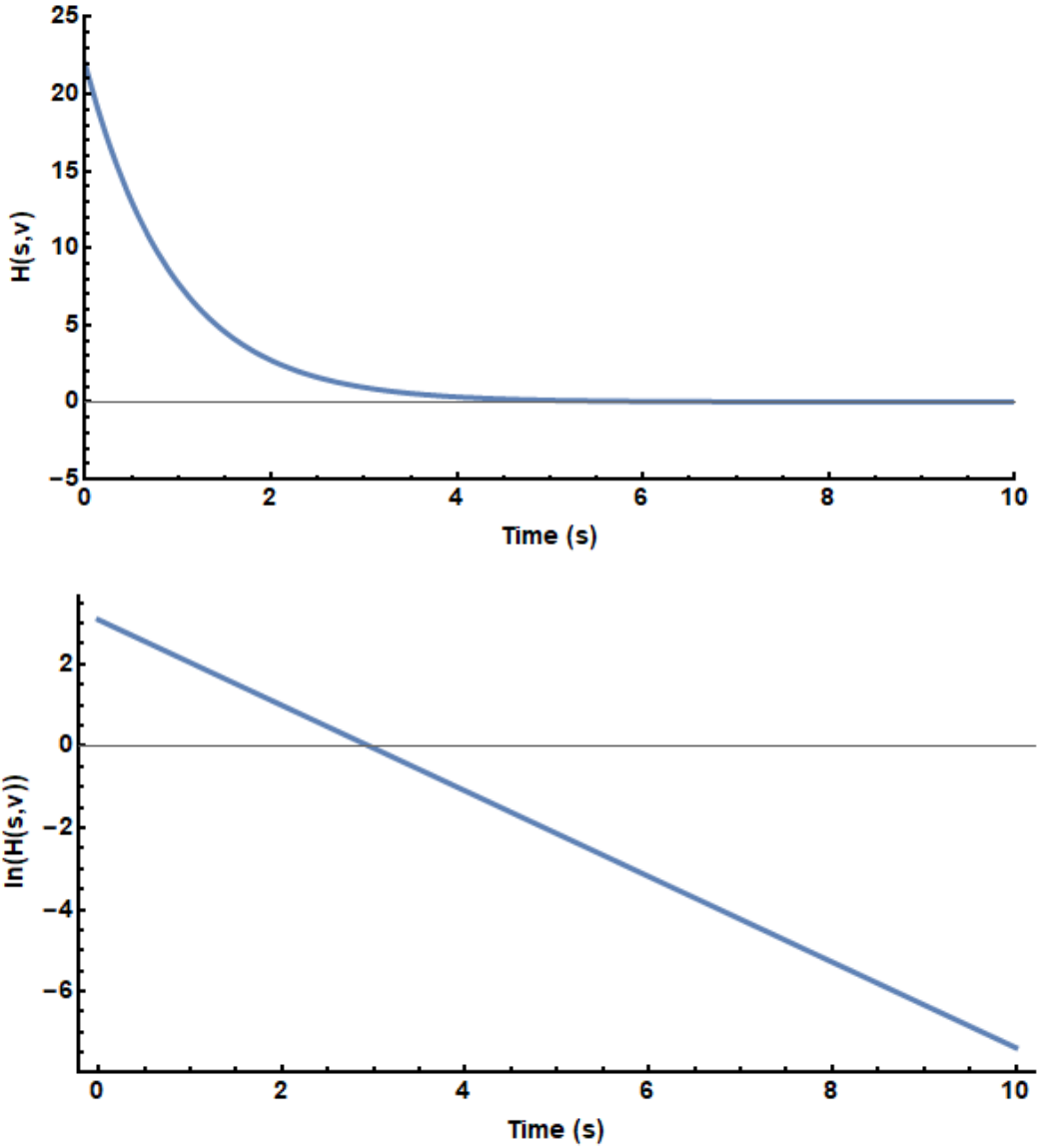

**Figure 4:** Convergence of the function $H(s,v)$ (top) and the graph of $\ln(H(s,v))$ (bottom) indicating exponential convergence.

**Example 3:** (Travelling waves of the macroscopic model (3.8)). To illustrate the results of Theorem 3.1, we consider a road with initial density and initial speed given by

$$\rho_0(x) = 0.1 + \begin{cases} 5x^2(x-1)^2, & x \in (0,1) \\ 0 & else \end{cases}$$

$$v_0(x) = 1 + \begin{cases} 8x^3(x-1)^3, & x \in (0,1) \\ 0 & else \end{cases}$$

These initial conditions indicate that there is a congestion belt on the interval $x \in (0,1)$ where vehicles are moving at lower speed and and accelerate again to a speed of $v^* = 1$ as the density decreases to a constant value. Values of $v^* = 1$ and $\omega = 1.2$ were used. In Figure 5 displays the speed profiles $v[t]$ (left) and the density profiles $\rho[t]$ (right) at different time instants $t = 0,1,...,5$.

This example illustrates the results of Theorem 3.1 which show that the speed converges exponentially to the speed set-point $v^*$; while the density converges to a travelling wave.

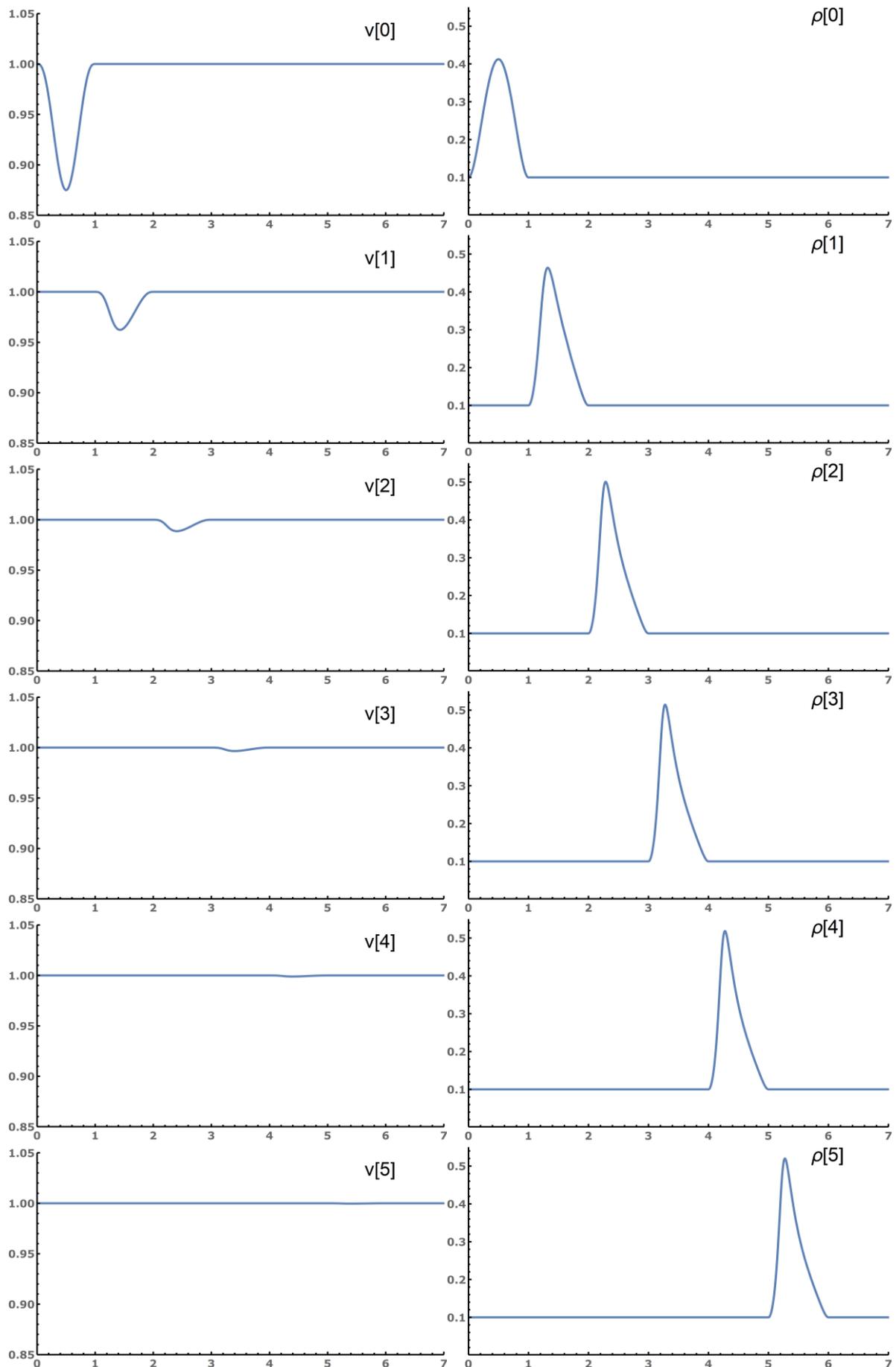

**Figure 5:** Exponential convergence of speed profiles $v[t]$ (left) and convergence of density profiles $\rho[t]$ to a travelling wave (right).

# 5. Numerical Investigation of Robustness

In this section we examine the sensitivity of system (2.6) to external disturbances. In particular, we consider how speed disturbances from the first vehicle propagate along a string of $n$ vehicles. We consider system (2.6) with the deviation of the speed of the leading vehicle from the equilibrium speed being the disturbance. The model becomes

$$\begin{aligned}
\dot{s}_2 &= v^* - v_2 + d \\
\dot{s}_i &= v_{i-1} - v_i, \quad i = 3,...,n \\
\dot{v}_i &= -k_i(s_i, s_{i+1})(v_i - v^*) + V'(s_i) - V'(s_{i+1}), \quad i = 2,...,n-1 \\
\dot{v}_n &= -k_n(s_n)(v_n - v^*) + V'(s_n)
\end{aligned} \quad (5.1)$$

where $d \in L^\infty(\Re_+)$ is given by $d = v_1 - v^*$, $v^* \in [0, v_{max}]$ is the desired speed set-point and $k_i, V$ are defined in (2.2), (2.3), respectively. We examine by simulations the attenuation of disturbances that satisfy

$$d(t) \in [-v^*, v_{max} - v^*], \ t \geq 0. \quad (5.2)$$

To quantify how the amplitudes of disturbances are attenuated by the bidirectional inviscid model (5.1), we define the *amplification factor* for the speed

$$\gamma_{i,n} = \frac{\|v_i - v^*\|_\infty}{\|d\|_\infty}, \quad i = 2,...,n, \ t \geq 0 \quad (5.3)$$

and *the amplification factor* for the spacing

$$\delta_{i,n} = \frac{\|V'(s_i)\|_\infty}{\|d\|_\infty}, \quad i = 2,...,n, \ t \geq 0. \quad (5.4)$$

We examine next the effects of sinusoidal disturbances

$$d(t) = \alpha \sin(\bar{\omega} t) \quad (5.5)$$

acting on system (5.1) for various frequencies $\bar{\omega} > 0$ and constant amplitude $\alpha \in [-v^*, v_{max} - v^*]$. For the following results, we have initialized the system at the equilibrium $v_i(0) = v^* = 30 m/s$, $s_i(0) = \lambda = 61m$, $i = 2,...,n$, and selected $V$ as in (4.1) with $v_{max} = 35 m/s$, $\mu = 2$, $L = 5.1m$, and $f$ as in (4.2) with $\varepsilon = 2$. We also consider the attenuation of disturbances of the nonlinear adaptive cruise controller with the Follow-the-Leader (FtL) model from [22]:

$$\begin{aligned}
\dot{s}_2 &= v^* - v_2 + d \\
\dot{v}_2 &= (k - \bar{g}(s_2))G(s_2) + \bar{g}(s_2)(v^* + d) - kv_2 \\
\dot{s}_i &= v_{i-1} - v_i \\
\dot{v}_i &= (k - \bar{g}(s_i))G(s_i) + \bar{g}(s_i)v_{i-1} - kv_i
\end{aligned}, \ i = 3,...,n \quad (5.6)$$

where $G(s) := \int_a^s \bar{g}(l)dl$, $a > 0$, and $\bar{g}(r)$ defined by

$$\bar{g}(r) = \begin{cases} 0 & r \leq b \\ (r-b) & \beta < r \leq g_{max} + b \\ g_{max} & g_{max} + b < r \leq \zeta \\ g_{max} \exp(\zeta - r) & r > \zeta \end{cases}$$

with $\beta > a > 0$ and $k > g_{max} > 0$. It is known that model (5.6) presents strong disturbance attenuation properties because it is string stable both in the $L^2$ norm and the $L^\infty$ norm (see [22]) and therefore can be used as a basis for comparison. Values of $a = 5.1$, $k = 1.2$, $\beta = 34.4$, $\zeta = 64.43$, and $g_{max} = 1.15$ were used. In the following simulation examples, both models (5.1) and (5.6) are initialized at the same equilibrium position.

**Example 4**: In this simulation example, we study the attenuation of disturbances for various frequencies $\bar{\omega}$ along a string of $n = 5$ vehicles, for the models (5.1) and the FtL model (5.6). We select $\alpha = -2.5$ and consider several frequencies $\bar{\omega}$. Figure 6 displays, with solid lines, the speed amplification factors $\gamma_{i,5}$, $i = 2,...,5$, for the inviscid model (5.1); and, with dashed lines, the values of $\gamma_{i,5}$, $i = 2,...,5$, for the model (5.6). It is seen that the disturbances dissipate faster for model (5.1), as the frequency of the disturbance increases. Figure 7 shows that the spacing factors $\delta_{i,5}$, $i = 2,...,5$, attain larger values for the inviscid model (5.1) compared to the FtL model. On the other hand, disturbances attenuate faster along the string of vehicles for the model (5.1). For clarity of presentation, only the factors $\delta_{i,5}$, $i = 3,4$ are shown.

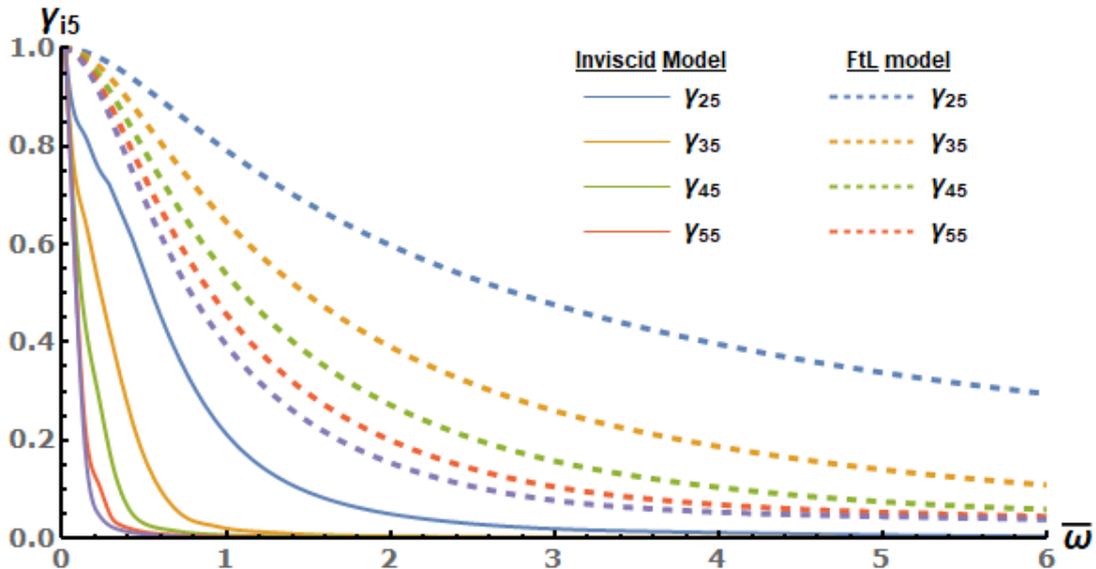

**Figure 6:** The effects of disturbance on speed $\gamma_{i,5}$ for a string of 5 vehicles and for several frequencies $\bar{\omega}$ for the inviscid model (5.1) (solid line) and the FtL model (5.6) (dashed line).

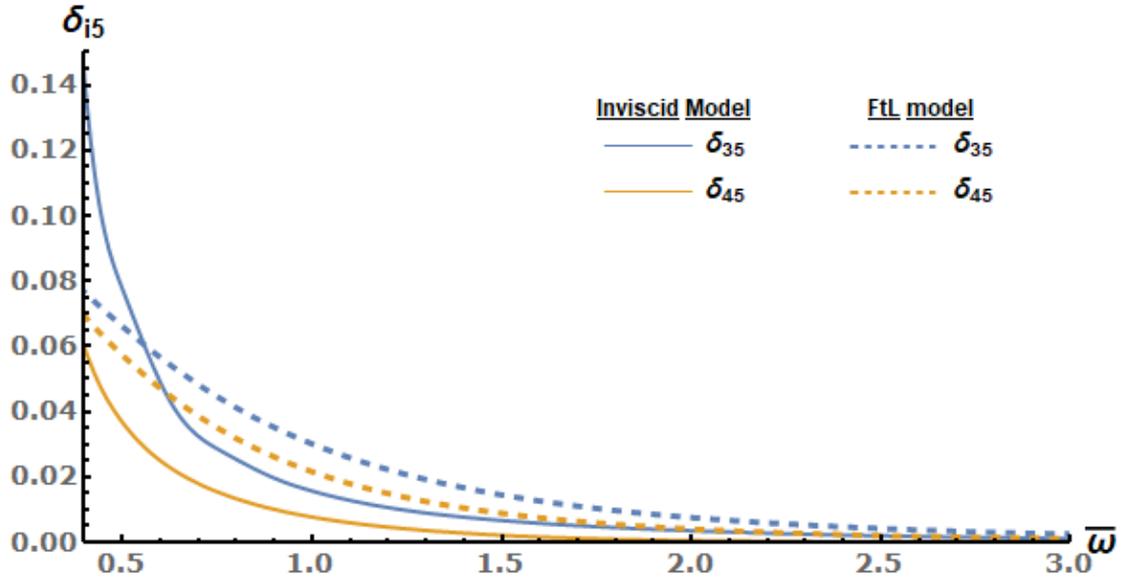

**Figure 7:** The effects of disturbance on speed $\delta_{i,5}$ for a string of 5 vehicles and for several frequences $\bar{\omega}$ for the inviscid model (5.1) (solid lines) and the FtL model (5.6) (dashed lines).

**Example 5:** We examine now, how disturbances with frequency $\bar{\omega}=0.25$ and amplitude $\alpha=-2.5$, dissipate along a string of $n=16$ vehicles. The speed amplification factors $\gamma_{i,16}$, $i=2,...,16$, are shown in Figure 8, which demonstrates that fluctuations on speed dissipate as they propagate backward along the string of vehicles. Finally, Figure 8 shows the spacing amplification factors $\delta_{i,n}$ which, while attaining larger values for the first 6 vehicles, dissipate faster as the number of vehicles increases.

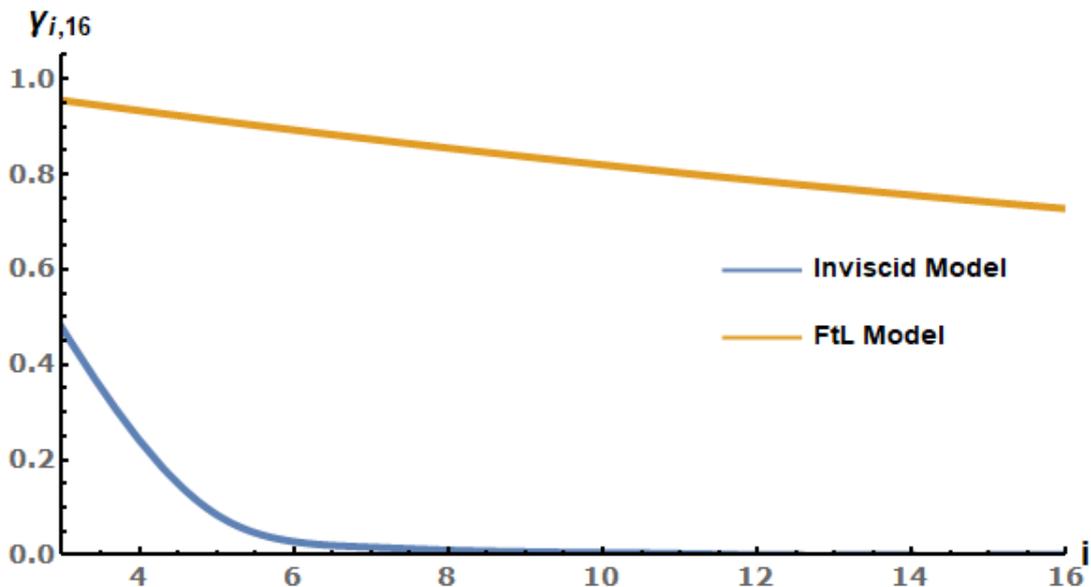

**Figure 8:** Amplification factors $\gamma_{i,16}$, $i=2,...,16$ for the inviscid model (5.1) and for the FtL model (5.6) for $\bar{\omega}=0.25$.

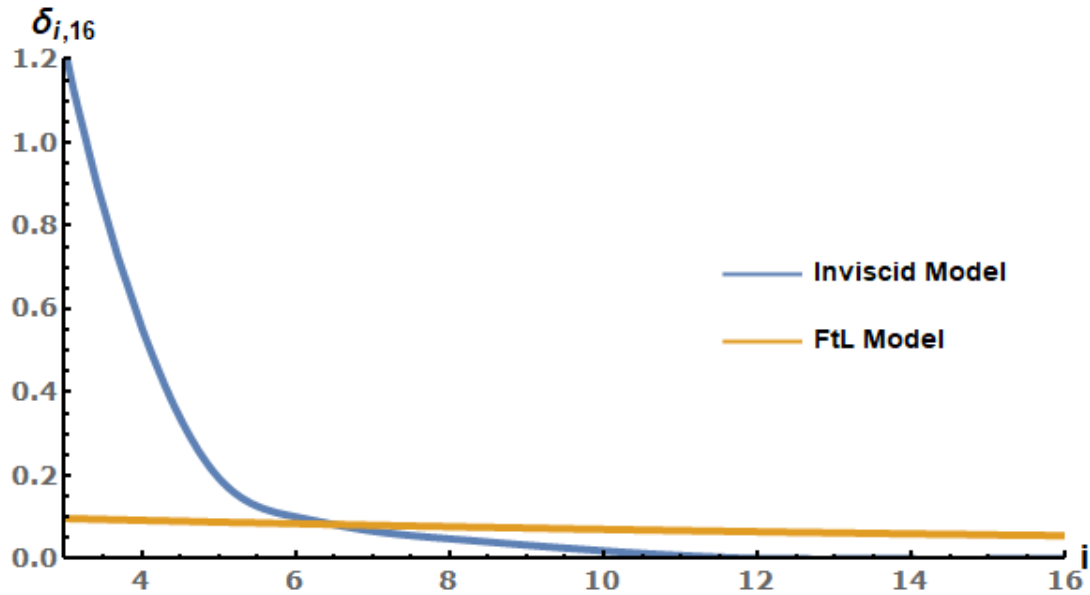

**Figure 9:** Amplification factors $\delta_{i,16}$, $i = 2,...,16$ for the inviscid model (5.1) and for the FtL model (5.6) for $\bar{\omega} = 0.25$.

The dissipation of the disturbances along the string of vehicle is even more evident by observing the amplification factors $\gamma_{n,n}$ and $\delta_{n,n}$ of the last vehicle in a string $n$ of vehicles for increasing values of $n$. Figure 10 and Figure 11 demonstrate the attentuation along the string of vehicles, as the total number of vehicles increases. It can be clearly seen that, for the inviscid model (5.1), disturbances dissipate much faster than the FtL model (5.6). In particular, for the cases of $n = 20$ and $n = 25$, both $\gamma_{n,n} = 0$ and $\delta_{n,n} = 0$ for the bidirectional model (5.1).

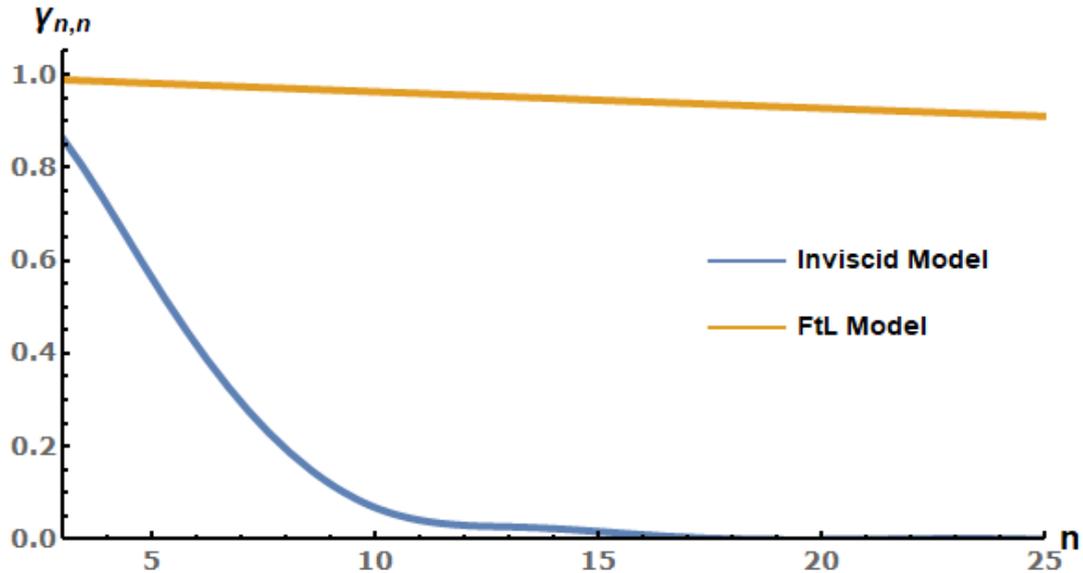

**Figure 10:** The effects of disturbances on the last vehicle on a string of $n$ vehicles for the inviscid model (5.1) and the FtL model (5.6) for $\bar{\omega} = 0.1$.

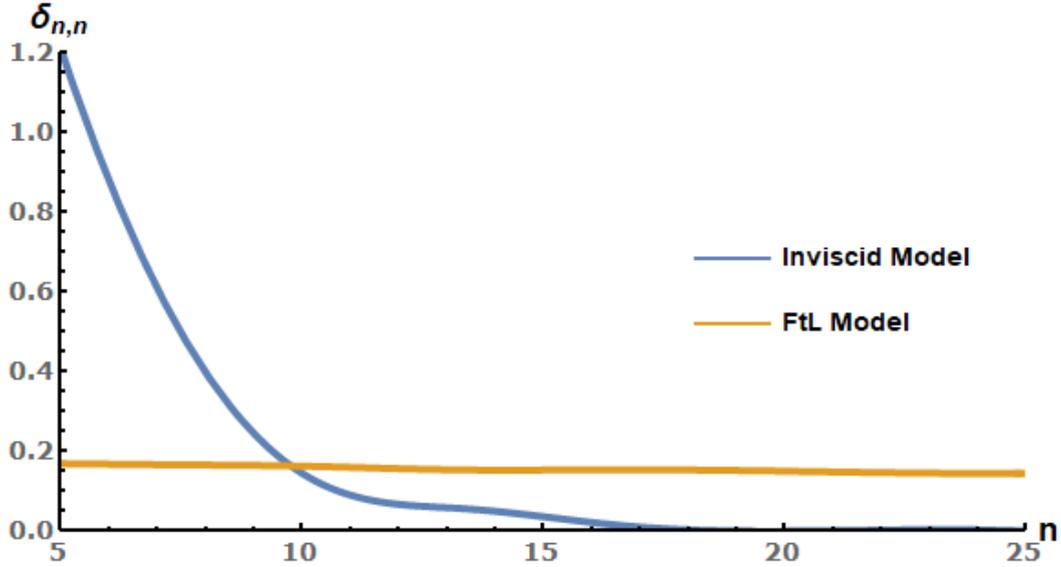

**Figure 11:** The effects of disturbances on the last vehicle on a string of $n$ vehicles for the inviscid model (5.1) and the FtL model (5.6) for $\bar{\omega}=0.1$.

From the previous simulations, we have seen that the factors $\gamma_{i,n}$, $i=2,...,n$, are much smaller for the inviscid model (5.1) than for the FtL model (5.6), and decrease as the frequency $\bar{\omega}$ and the number of vehicles $n$ increase. Finally, we have observed that the factors $\delta_{i,n}$, $i=2,...,n$, are larger for the inviscid model, than the FtL model (5.6) for small $\bar{\omega}$, but dissipate at a much higher rate as $\bar{\omega}$ and $n$ are increased.

## 6. Proofs of Main Results

We first provide a technical result that is invoked many times for the proofs of the main results.

**Lemma 6.1:** *Let $\lambda > L > 0$ and define $H:\Omega \to \mathfrak{R}_+$ by means of (2.10) where $\Omega$ is given by (2.7). Then, there exists a non-increasing function $\rho:\mathfrak{R}_+ \to (L,\lambda]$ that the following inequalities hold for all $(s,v) \in \Omega$:*

$$s_i \geq \rho(H(s,v)), \text{ for all } i = 2,...,n \qquad (6.1)$$

**Proof of Lemma 6.1:** Notice first that, due to (2.3), $V(r)$ is decreasing for $r \in (L,\lambda)$ with $V(r)=0$ for $r \geq \lambda$. Thus, the inverse function $V^{-1}:\mathfrak{R}_+ \to (L,\lambda]$ of $V$ (restricted on the interval $(L,\lambda]$), exists, and due to continuity of $V$, is continuous, decreasing and satisfies $V^{-1}(0) = \lambda$ and $\lim_{l \to +\infty}\left(V^{-1}(l)\right) = L$. Next, using (2.10), we get that $V(s_i) \leq H(s,v)$ for all $(s,v) \in \Omega$ and $i = 2,...,n$. It follows from the previous inequality that $s_i \geq V^{-1}(H(s,v)) > L$. The previous inequality shows that (6.1) holds with $\rho(r):=V^{-1}(r)$ for $r \geq 0$. This concludes the proof. ◁

**Proof of Theorem 2.1:** First, we show that the solution of the model is defined for all $t \geq 0$. Let $(s_0,v_0) = (s_{2,0},...,s_{n,0},v_{1,0},...,v_{n,0}) \in \Omega$ be given. Let $\varepsilon > 0$ and consider the open set

$$\Omega_\varepsilon := \left\{ (s_2,...,s_n,v_1,...,v_n) \in \mathfrak{R}^{2n-1} : \min_{i=2,...,n}(s_i) > L, \max_{i=1,...,n}(v_i) < v_{max} + \varepsilon, \min_{i=1,...,n}(v_i) > -\varepsilon \right\}, \quad (6.2)$$

and notice that $\Omega \subset \Omega_\varepsilon$. Due to the fact that $k_i$ (defined by (2.2)) and $V'$ are locally Lipschitz, there exists $t_{max} \in (0,+\infty]$ such that the unique solution $(s(t),v(t)) = (s_2(t),...,s_n(t),v_1(t),...,v_n(t))$ of the initial value problem (2.6) with initial conditions $(s(0),v(0)) = (s_0,v_0) \in \Omega \subset \Omega_\varepsilon$ is defined on $[0,t_{max})$ and satisfies $(s(t),v(t)) \in \Omega_\varepsilon$ for all $t \in [0,t_{max})$. If $t_{max} < +\infty$, then, we necessarily have that either $\limsup_{t \to t_{max}^-}(|(s(t),v(t))|) = +\infty$ or $\lim_{t \to t_{max}^-}(dist((s(t),v(t)),\partial \Omega_\varepsilon)) = 0$.

Notice first that, since $(s(t),v(t)) \in \Omega_\varepsilon$ for all $t \in [0,t_{max})$, it follows from (2.6) and definitions (2.7), (6.2) that $\dot s_i \leq v_{max} + 2\varepsilon$ for all $t \in [0,t_{max})$ and $i = 2,...,n$ which in turn implies that

$$s_i(t) \leq s_i(0) + (v_{max} + 2\varepsilon)t \quad (6.3)$$

for all $t \in [0,t_{max})$ and $i = 2,...,n$. Moreover, for all $(s,v) \in \Omega_\varepsilon$, inequality (2.11) holds and implies that

$$H(s(t),v(t)) \leq H(s(0),v(0)) \quad \text{for all } t \in [0,t_{max}). \quad (6.4)$$

In addition, due to Lemma 6.1, there exists a non-increasing function $\rho : \mathfrak{R}_+ \to (L,\lambda]$ such that inequality (6.1) holds for all $(s,v) \in \Omega$. It follows then, by (6.1) and (6.4) that for all $t \in [0,t_{max})$ and $i = 2,...,n$ the following holds

$$s_i(t) \geq \rho(H(s(t),v(t))) \geq \rho(H(s(0),v(0))) > L. \quad (6.5)$$

Define

$$b_1(s) := \max\{|V'(d)| : s \leq d \leq \lambda\} \text{ for } s \in (L,\lambda]$$

Then, for all $(s,v) \in \Omega$, we obtain from inequality (6.1) and the above definition that

$$|V'(s_i) - V'(s_{i+1})| \leq |V'(s_i)| + |V'(s_{i+1})| \leq 2b_1(\rho(H(s,v)))$$

The above inequality and definition (2.2) imply that

$$\mu \leq k_i(s_i,s_{i+1}) \leq \varphi(H(s,v)) \quad (6.6)$$

where $\varphi(r) := \mu + \max\{g(z) : |z| \leq 2b_1(\rho(s))\}$. Notice also that definitions (2.2), (2.4) and (2.5) imply that:

$$k_i(s_i)v^* \geq V'(s_{i+1}) - V'(s_i) \geq -k_i(s_i,s_{i+1})(v_{max} - v^*), \quad i = 2,...,n-1$$
$$k_1(s_2)v^* \geq -V'(s_2) \geq -k_1(s_2)(v_{max} - v^*) \quad (6.7)$$
$$k_n(s_n)v^* \geq V'(s_n) \geq -k_n(s_n)(v_{max} - v^*)$$

By using inequalities (6.7), (2.6) and the fact that $v_i \in [0,v_{max}]$, we get for all $(s,v) \in \Omega$

$$\begin{aligned}
&\dot{v}_i \geq -k_i(s_i, s_{i+1})(v_i - v^*) - k_i(s_i, s_{i+1})v^* = -k_i(s_i, s_{i+1})v_i \geq -k_i(s_i, s_{i+1})v_{\max}, \ i = 2, \ldots, n-1 \\
&\dot{v}_1 \geq -k_1(s_2)v_1 \geq -k_1(s_2)v_{\max}, \quad \dot{v}_n \geq -k_n(s_n)v_{\max} \\
&\dot{v}_i \leq -k_i(s_i, s_{i+1})(v_i - v^*) + k_i(s_i, s_{i+1})(v_i - v^*) = k_i(s_i, s_{i+1})(v_{\max} - v_i) \leq k_i(s_i, s_{i+1})v_{\max}, \ i = 2, \ldots, n-1 \\
&\dot{v}_1 \leq k_1(s_2)v_{\max}, \quad \dot{v}_n \leq k_n(s_n)v_{\max}
\end{aligned} \quad (6.8)$$

Then, due (2.6), (6.6), (6.7), (6.8), and (6.4) the following inequalities hold for all $t \in [0, t_{\max})$:

$$k_i(t)(v_{\max} - v_i(t)) \geq \dot{v}_i(t) \geq -k_i(t)v_i(t), \tag{6.9}$$

$$k_i(t) \leq M := \varphi(H(s(0), v(0))). \tag{6.10}$$

An immediate consequence of inequalities (6.9) and (6.10), the fact that $(s_0, v_0) \in \Omega$, definition (2.7) and the comparison principle in [23] is the following estimate

$$v_{\max} \geq v_i(0)\exp(-Mt) + (1 - \exp(-Mt))v_{\max} \geq v_i(t) \geq v_i(0)\exp(-Mt) \geq 0 \tag{6.11}$$

for all $t \in [0, t_{\max})$ and $i = 1, \ldots, n$.

Inequalities (6.5) and (6.11) imply that $(s(t), v(t)) \in \Omega$ for all $t \in [0, t_{\max})$ and in addition $t_{\max} = +\infty$. Indeed, we conclude that $t_{\max} = +\infty$, since inequalities (6.3), (6.5), and (6.11) imply that neither $\limsup_{t \to t_{\max}^-}(|(s(t), v(t))|) = +\infty$ nor $\lim_{t \to t_{\max}^-}(dist((s(t), v(t)), \partial \Omega_\varepsilon)) = 0$ can hold when $t_{\max} < +\infty$.

In order to prove inequality (2.12), we proceed as follows. Suppose that there exists $i \in \{2, \ldots, n\}$ and $T \geq 0$ for which $s_i(T) > \max(\lambda, s_i(0)) + \mu^{-1}v_{\max}$. Clearly, $T > 0$. Define:

$$A = \{t \in [0, T] : s_i(t) \leq \lambda\} \tag{6.12}$$

Suppose that $A \neq \emptyset$. Then we define $t^* = \sup(A)$ and clearly $t^* \in [0, T]$. Moreover, continuity of the mapping $t \to s_i(t)$ and the fact that $s_i(T) > \max(\lambda, s_i(0)) + \mu^{-1}v_{\max}$ imply that $t^* < T$. Furthermore, by definition of the set $A$ we have $s_i(t) \geq \lambda$ for all $t \in [t^*, T]$ and $s_i(t^*) = \lambda$. If $i \in \{2, \ldots, n-1\}$ then (recall (2.3) which implies that $V'(s_i(t)) = 0$ for all $t \in [t^*, T]$)

$$\dot{v}_i(t) = -k_i(t)\left(v_i(t) - v^*\right) - V'(s_{i+1}(t)), \text{ for all } t \in [t^*, T]$$

and consequently (due to (2.3)):

$$\dot{v}_i(t) \geq -k_i(t)\left(v_i(t) - v^*\right), \text{ for all } t \in [t^*, T].$$

The above differential inequality holds even if $i = n$ and implies that (using the comparison principle in [23]):

$$v_i(t) \geq v^* + \left(v_i(t^*) - v^*\right)\exp\left(-\int_{t^*}^{t} k_i(l)dl\right), \text{ for all } t \in [t^*, T]. \tag{6.13}$$

Moreover, we get from (2.6), (2.3) and the fact that $V'(s_i(t)) = 0$ for all $t \in [t^*, T]$:

$$\dot{v}_{i-1}(t) \leq -k_{i-1}(t)\left(v_{i-1}(t) - v^*\right), \text{ for all } t \in [t^*, T].$$

The above differential inequality implies that (using the comparison principle in [23]):

$$v_{i-1}(t) \leq v^* + \left(v_{i-1}(t^*) - v^*\right)\exp\left(-\int_{t^*}^{t} k_{i-1}(l)dl\right), \text{ for all } t \in [t^*, T]. \tag{6.14}$$

Thus, from (2.6), (6.13), and (6.14) we obtain for all $t \in [t^*, T]$

$$\dot{s}_i(t) \leq \left(v_{i-1}(t^*) - v^*\right)\exp\left(-\int_{t^*}^{t} k_{i-1}(l)dl\right) + \left(v^* - v_i(t^*)\right)\exp\left(-\int_{t^*}^{t} k_i(l)dl\right), \text{ for all } t \in [t^*, T]$$

Therefore, since $s_i(t^*) = \lambda$ we get for all $t \in [t^*, T]$:

$$s_i(t) \leq \lambda + \left(v_{i-1}(t^*) - v^*\right)\int_{t^*}^{t}\exp\left(-\int_{t^*}^{\tau} k_{i-1}(l)dl\right)d\tau + \left(v^* - v_i(t^*)\right)\int_{t^*}^{t}\exp\left(-\int_{t^*}^{\tau} k_i(l)dl\right)d\tau$$

$$\leq \lambda + \max\left(0, v_{i-1}(t^*) - v^*\right)\int_{t^*}^{t}\exp\left(-\int_{t^*}^{\tau} k_{i-1}(l)dl\right)d\tau + \max\left(0, v^* - v_i(t^*)\right)\int_{t^*}^{t}\exp\left(-\int_{t^*}^{\tau} k_i(l)dl\right)d\tau$$

Since $k_i(t) \geq \mu$, $k_{i-1}(t) \geq \mu$ (a consequence of (2.2), (2.4), (2.5)) we also have for all $t \in [t^*, T]$

$$s_i(t) \leq \lambda + \left(\max\left(0, v_{i-1}(t^*) - v^*\right) + \max\left(0, v^* - v_i(t^*)\right)\right)\int_{t^*}^{t}\exp(-\mu(\tau - t^*))d\tau$$

$$\leq \lambda + \mu^{-1}\left(\max\left(0, v_{i-1}(t^*) - v^*\right) + \max\left(0, v^* - v_i(t^*)\right)\right) \tag{6.15}$$

Since $\max\left(0, v_{i-1}(t^*) - v^*\right) + \max\left(0, v^* - v_i(t^*)\right) \leq v_{\max}$ we obtain from (6.15) that $s_i(T) \leq \lambda + \mu^{-1}v_{\max}$. This contradicts the assumption $s_i(T) > \max(\lambda, s_i(0)) + \mu^{-1}v_{\max}$.

If $A = \emptyset$ then we perform exactly the above steps with $t^*$ replaced by 0 and we obtain the estimate $s_i(T) \leq s_i(0) + \mu^{-1}v_{\max}$, which contradicts the assumption $s_i(T) > \max(\lambda, s_i(0)) + \mu^{-1}v_{\max}$.

Finally, we conclude the proof of Theorem 2.1 by showing that $\lim_{t \to +\infty}(v_i(t)) = v^*$ for all $i = 1, ..., n$ and $\lim_{t \to +\infty}(V(s_i(t))) = 0$ for all $i = 2, ..., n$. To that end, let $(s_0, v_0) \in \Omega$, and define the set

$$\Omega_c := \overline{\bigcup_{t \geq 0}\{(s(t), v(t))\}} \tag{6.16}$$

and notice that since for all $i = 2,...,n$, $s_i(t)$ satisfies (2.12), (6.5) for all $t \geq 0$ and since $v_i(t) \in [0, v_{max}]$ for all $t \geq 0$, $i = 1,...,n$, the set $\Omega_c$ is compact, positively invariant and satisfies $\Omega_c \subset \Omega$. Define also the set

$$Q := \{(s,v) \in \Omega_c : \dot{H}(s,v) = 0\} = \{(s,v) \in \Omega_c : v_i = v^*, i = 1,...,n\} \qquad (6.17)$$

Then, due to LaSalle's Invariance Principle, see [23], the state $(s(t), v(t))$ approaches the largest invariant set in $Q$ as $t \to +\infty$. It follows from properties (2.3), the model (2.6) and (6.17) that the largest invariant set contained in $Q$ is $M = \Omega_c \cap S$. Hence, $\lim_{t \to +\infty} (v_i(t)) = v^*$ for all $i = 1,...,n$ and (due to (2.3)) $\lim_{t \to +\infty} (V(s_i(t))) = 0$ for all $i = 2,...,n$. The proof is complete. ◁

**Proof of Theorem 2.2:** Let $\beta > 0$ be given (arbitrary). The first part of the proof is devoted to the construction of the function $R \in C^1(\mathfrak{R}_+; (0, +\infty))$. We proceed by proving some useful claims.

**Claim 1:** *The following inequalities hold for all $(s,v) \in \Omega$:*

$$k_1(s_2) \leq \varphi(H(s,v)), \ k_i(s_i, s_{i+1}) \leq \varphi(H(s,v)), \text{ for } i = 2,...,n-1, \ k_n(s_n) \leq \varphi(H(s,v)) \qquad (6.18)$$

*where $\varphi : \mathfrak{R}_+ \to \mathfrak{R}_+$ is the non-decreasing function defined by*

$$\varphi(r) := \mu + \max\{g(z) : |z| \leq 2b_1(\rho(r))\} \qquad (6.19)$$

*with*

$$b_1(r) := \max\{|V'(d)| : r \leq d \leq \lambda\} \text{ for } r \in (L, \lambda] \qquad (6.20)$$

Claim 1 is a direct consequence of inequality (6.1), the definitions (2.2), (2.4) of $k_i$ ($i = 1,...,n$) and the fact that $V'(d) = 0$ for $d \geq \lambda$.

**Claim 2:** *The following inequalities hold for all $(s,v) \in \Omega$:*

$$V''(s_i) \leq \tilde{\varphi}(H(s,v)), \text{ for } i = 2,...,n \qquad (6.21)$$

*where $\tilde{\varphi} : \mathfrak{R}_+ \to \mathfrak{R}_+$ is the non-decreasing function defined by*

$$\tilde{\varphi}(r) := b_2(\rho(r)) \qquad (6.22)$$

*with*

$$b_2(r) := \max\{V''(d) : r \leq d \leq \lambda\}, \text{ for } r \in (L, \lambda] \qquad (6.23)$$

Claim 2 is a direct consequence of inequality (6.1) and the fact that $V''(d) = 0$ for $d \geq \lambda$.

**Claim 3:** *There exists a non-decreasing function $\gamma : \mathfrak{R}_+ \to \mathfrak{R}_+$ such that the following inequality holds for all $(s,v) \in \Omega$:*

$$(V'(s_i))^2 \leq \gamma(H(s,v))V(s_i), \text{ for } i = 2,...,n \qquad (6.24)$$

**Proof of Claim 3:** Since $\lim_{d \to \lambda^-} \left( \frac{(V'(d))^2}{V(d)} \right) = 2 \lim_{d \to \lambda^-} (V''(d)) = 0$, it follows that the function $h:(L,+\infty) \to \Re_+$ defined by

$$h(r) := \sup_{r \leq d < \lambda} \left( \frac{(V'(d))^2}{V(d)} \right), \text{ for } r \in (L,\lambda) \text{ and } h(r) := 0 \text{ for } r \geq \lambda \tag{6.25}$$

is a well-defined, non-increasing function which satisfies the following inequality for all $d \in (L,\lambda)$:

$$(V'(d))^2 \leq h(d)V(d) \tag{6.26}$$

Since $V(d) = V'(d) = V''(d) = 0$ for $d \geq \lambda$, it follows that (6.26) holds for all $d > L$. Inequality (6.24) with $\gamma(r) := h(\rho(r))$ is a direct consequence of (6.26) and (6.1). ◁

Define

$$\tilde{R}(r) := 2 + \frac{4^{2n}}{2}\gamma(r) + 4^n \left( \varphi(r) + \frac{7}{2\mu}\tilde{\varphi}(r) \right) + \frac{\beta}{\mu}, \text{ for } r \geq 0 \tag{6.27}$$

Clearly, the function $\tilde{R}: \Re_+ \to \Re_+$ defined by (6.27) is a non-decreasing function. Consequently, there exists a non-decreasing function $R \in C^1(\Re_+;(0,+\infty))$ that satisfies

$$R(r) \geq \tilde{R}(r), \text{ for } r \geq 0 \tag{6.28}$$

Having constructed $R \in C^1(\Re_+;(0,+\infty))$, we proceed by showing the validity of inequalities (2.13), (2.14).

Using (6.27), (6.28), (6.18), (6.21), (6.24), it follows that the following inequalities hold for all $(s,v) \in \Omega$:

$$R(H(s,v)) \geq 4^i \left( k_i(s_i, s_{i+1}) + \frac{3}{2\mu}V''(s_i) + \frac{2}{\mu}V''(s_{i+1}) \right) + \frac{\beta}{\mu}, \quad i = 2,...,n-1$$

$$R(H(s,v)) \geq \frac{8}{\mu}V''(s_2) + \frac{\beta}{\mu} \tag{6.29}$$

$$R(H(s,v)) \geq 4^n \left( k_n(s_n) + \frac{3}{2\mu}V''(s_n) \right) + \frac{\beta}{\mu}$$

$$R(H(s,v)) \geq 2 + \frac{4^{2n}}{2}\gamma(H(s,v)) \tag{6.30}$$

In what follows we omit the arguments of the functions $k_i, i = 1,...,n$ (for simplicity). Definition (2.15) and the equations of the model give for all $(s,v) \in \Omega$:

$$\dot{W} = -\left(R'(H(s,v))H(s,v) + R(H(s,v))\right)\sum_{i=1}^{n}k_i\left(v_i - v^*\right)^2 - \sum_{i=2}^{n}4^i\left(V'(s_i)\right)^2$$

$$+\sum_{i=2}^{n-1}4^i k_i\left(v_i - v^*\right)V'(s_i) + \sum_{i=2}^{n-1}4^i V'(s_i)V'(s_{i+1}) - \sum_{i=2}^{n}4^i V''(s_i)\left(v_{i-1} - v^*\right)\left(v_i - v^*\right) \quad (6.31)$$

$$+\sum_{i=2}^{n}4^i V''(s_i)\left(v_i - v^*\right)^2 + 4^n k_n\left(v_n - v^*\right)V'(s_n)$$

Since $R'(r) \geq 0$ for all $r \geq 0$ (recall that $R$ is non-decreasing) and using the inequalities

$$k_i\left(v_i - v^*\right)V'(s_i) \leq \frac{1}{4}\left(V'(s_i)\right)^2 + k_i^2\left(v_i - v^*\right)^2$$

$$\left|v_{i-1} - v^*\right|\left|v_i - v^*\right| \leq \frac{1}{2}\left(v_{i-1} - v^*\right)^2 + \frac{1}{2}\left(v_i - v^*\right)^2$$

$$V'(s_i)V'(s_{i+1}) \leq \frac{1}{2}\left(V'(s_i)\right)^2 + \frac{1}{2}\left(V'(s_{i+1})\right)^2$$

we get from (6.31) the following inequality for all $(s,v) \in \Omega$:

$$\dot{W} \leq -\left(R(H(s,v)) - \frac{8V''(s_2)}{k_1}\right)k_1\left(v_1 - v^*\right)^2$$

$$-\sum_{i=2}^{n-1}k_i\left(R(H(s,v)) - 4^i k_i - \frac{3}{2k_i}4^i V''(s_i) - \frac{1}{2k_i}4^{i+1}V''(s_{i+1})\right)\left(v_i - v^*\right)^2 \quad (6.32)$$

$$-k_n\left(R(H(s,v)) - 4^n k_n - \frac{3}{2k_n}4^n V''(s_n)\right)\left(v_n - v^*\right)^2$$

$$-4\left(V'(s_2)\right)^2 - \frac{1}{8}\sum_{i=3}^{n-1}4^i\left(V'(s_i)\right)^2 - 5\frac{4^n}{8}\left(V'(s_n)\right)^2$$

Inequality (6.32) in conjunction with inequalities (6.29) and the fact that $k_i \geq \mu$ for $i = 1,\ldots,n$, give inequality (2.14).

We also obtain from definition (2.15) by completing the squares and using (6.24):

$$W(s,v) \geq R(H(s,v))H(s,v) - \frac{1}{2}\sum_{i=2}^{n}\left(v_i - v^*\right)^2 - \frac{1}{2}\sum_{i=2}^{n}4^{2i}\left(V'(s_i)\right)^2$$

$$\geq R(H(s,v))H(s,v) - \frac{1}{2}\sum_{i=2}^{n}\left(v_i - v^*\right)^2 - \frac{4^{2n}}{2}\sum_{i=2}^{n}\left(V'(s_i)\right)^2$$

$$\geq R(H(s,v))H(s,v) - \frac{1}{2}\sum_{i=2}^{n}\left(v_i - v^*\right)^2 - \frac{4^{2n}}{2}\gamma(H(s,v))\sum_{i=2}^{n}V(s_i)$$

$$\geq R(H(s,v))H(s,v) - \left(1 + \frac{4^{2n}}{2}\gamma(H(s,v))\right)\left(\frac{1}{2}\sum_{i=2}^{n}\left(v_i - v^*\right)^2 + \sum_{i=2}^{n}V(s_i)\right)$$

$$= \left(R(H(s,v)) - 1 - \frac{4^{2n}}{2}\gamma(H(s,v))\right)H(s,v)$$

The above inequality in conjunction with (6.30) gives the left inequality (2.13).

Finally, we obtain from definition (2.15) by completing the squares and using (6.24):

$$\begin{aligned} W(s,v) &\leq R(H(s,v))H(s,v) + \frac{1}{2}\sum_{i=2}^{n}(v_i - v^*)^2 + \frac{1}{2}\sum_{i=2}^{n} 4^{2i}(V'(s_i))^2 \\ &\leq R(H(s,v))H(s,v) + \frac{1}{2}\sum_{i=2}^{n}(v_i - v^*)^2 + \frac{4^{2n}}{2}\sum_{i=2}^{n}(V'(s_i))^2 \\ &\leq R(H(s,v))H(s,v) + \frac{1}{2}\sum_{i=2}^{n}(v_i - v^*)^2 + \frac{4^{2n}}{2}\gamma(H(s,v))\sum_{i=2}^{n}V(s_i) \\ &\leq R(H(s,v))H(s,v) + \left(1 + \frac{4^{2n}}{2}\gamma(H(s,v))\right)\left(\frac{1}{2}\sum_{i=2}^{n}(v_i - v^*)^2 + \sum_{i=2}^{n}V(s_i)\right) \\ &= \left(1 + R(H(s,v)) + \frac{4^{2n}}{2}\gamma(H(s,v))\right)H(s,v) \end{aligned}$$

The above inequality gives the right inequality (2.13) for any non-decreasing function $\kappa \in C^0(\Re_+;(0,+\infty))$ that satisfies $\kappa(r) \geq 1 + R(r) + \frac{4^{2n}}{2}\gamma(r)$ for all $r \geq 0$ (e.g. the function $\kappa(r) = 1 + R(r) + \frac{4^{2n}}{2}\int_r^{r+1}\gamma(l)dl$). The proof is complete. ◁

**Proof of Theorem 2.3:** Define the operator $P:\Re^{n-1} \to \Re^{n-1}$ defined by

$$\Re^{n-1} \ni (s_2,...,s_n) \to y = (y_2,...,y_n) = P(s) \text{ with } y_i = \begin{cases} s_i & \text{if } s_i \leq \lambda \\ \lambda & \text{if } s_i > \lambda \end{cases} \quad (6.33)$$

Due to the fact that (2.3) holds, we obtain from (2.6) and definitions (2.10), (2.15), (6.33) for all $(s,v) \in \Omega$:

$$\begin{aligned} \dot{W}(s,v) &= \dot{W}(P(s),v) \\ W(s,v) &= W(P(s),v) \end{aligned} \quad (6.34)$$

Define the parameterized family of sets with parameter $a \in (L,\lambda]$:

$$\Theta(a) = \left\{ (s_2,...,s_n,v_1,...,v_n) \in \Re^{2n-1} : \min_{i=2,...,n}(s_i) \geq a, \max_{i=1,...,n}(v_i) \leq v_{\max}, \min_{i=1,...,n}(v_i) \geq 0, \max_{i=2,...,n}(s_i) \leq \lambda \right\} \quad (6.35)$$

Notice that $\Theta(a) \subseteq \Omega$ for all $a \in (L,\lambda]$ and definition (6.35) in conjunction with inequalities (6.1) and (2.13) imply that

$$(P(s),v) \in \Theta(\rho(W(s,v))) \text{ for all } (s,v) \in \Omega \quad (6.36)$$

Moreover, notice that $\Theta(a) \subseteq \Omega$ is a compact set for all $a \in (L,\lambda]$.

Define the function $\tilde{\varphi}:\Re_+ \to \Re_+$ by

$$\tilde{\varphi}(r) := \inf\left\{-\dot{W}(s,v) : (s,v) \in \Omega, r \leq W(s,v) \leq 1\right\}, \text{ for } r \in [0,1] \quad (6.37)$$

$$\tilde{\varphi}(r) := \inf\left\{-\dot{W}(s,v) : (s,v) \in \Omega, 1 \leq W(s,v) \leq r\right\}, \text{ for } r > 1 \quad (6.38)$$

Clearly, inequality (2.14) implies that $\tilde{\varphi}(r) \geq 0$ for all $r \geq 0$ and $\tilde{\varphi}(0) = 0$. Furthermore, we obtain from (6.34) and (6.36) by virtue of continuity of $W(s,v)$ and $\dot{W}(s,v)$ and compactness of $\Theta(\rho(r)) \subseteq \Omega$ for all $r \geq 0$:

$$\begin{aligned}
\tilde{\varphi}(r) &= \inf\left\{-\dot{W}(s,v) : (s,v) \in \Omega, r \leq W(s,v) \leq 1\right\} \\
&= \inf\left\{-\dot{W}(P(s),v) : (s,v) \in \Omega, r \leq W(P(s),v) \leq 1\right\} \\
&\geq \inf\left\{-\dot{W}(y,v) : (y,v) \in \Theta(\rho(1)), r \leq W(y,v) \leq 1\right\} \\
&= \min\left\{-\dot{W}(y,v) : (y,v) \in \Theta(\rho(1)), r \leq W(y,v) \leq 1\right\}
\end{aligned}$$

for $r \in [0,1]$ and

$$\begin{aligned}
\tilde{\varphi}(r) &:= \inf\left\{-\dot{W}(s,v) : (s,v) \in \Omega, 1 \leq W(s,v) \leq r\right\} \\
&= \inf\left\{-\dot{W}(P(s),v) : (s,v) \in \Omega, 1 \leq W(P(s),v) \leq r\right\} \\
&\geq \inf\left\{-\dot{W}(y,v) : (s,v) \in \Theta(\rho(r)), 1 \leq W(y,v) \leq r\right\} \\
&= \min\left\{-\dot{W}(y,v) : (s,v) \in \Theta(\rho(r)), 1 \leq W(y,v) \leq r\right\}
\end{aligned}$$

for $r > 1$.

The above inequalities in conjunction with (2.14), continuity of $W(s,v)$ and $\dot{W}(s,v)$ and compactness of $\Theta(\rho(r)) \subseteq \Omega$ for all $r \geq 0$ imply that

$$\tilde{\varphi}(r) > 0 \quad \text{for all } r > 0 \quad (6.39)$$

Therefore, the function $\tilde{\varphi} : \Re_+ \to \Re_+$ defined by (6.37), (6.38) is positive definite. Notice that definitions (6.37), (6.38) also guarantee that $\tilde{\varphi} : \Re_+ \to \Re_+$ is non-decreasing on $[0,1]$ and non-increasing on $(1,+\infty)$. The following inequality is also a direct consequence of definitions (6.37), (6.38):

$$\dot{W}(s,v) \leq -\tilde{\varphi}(W(s,v)), \text{ for all } (s,v) \in \Omega \quad (6.40)$$

Following exactly the same methodology as in the proof of Proposition 2.2 on page 107 in [19], we can construct a locally Lipschitz, positive definite function $\varphi : \Re_+ \to \Re_+$ that satisfies $\tilde{\varphi}(r) \geq \varphi(r)$ for all $r \geq 0$. Lemma 2.13 on page 80 in [19] in conjunction with inequalities (2.14), (6.39), (6.40) and the fact that $\tilde{\varphi}(r) \geq \varphi(r)$ for all $r \geq 0$ implies the existence of $\sigma \in KL$ that satisfies the right inequality (2.16) for all $t \geq 0$.

In order to finish the proof, we need to establish the existence of a function $a \in K_\infty$ so that $a(dist((s,v),S)) \leq H(s,v)$ for all $(s,v) \in \Omega$. Indeed, in this case the use of (2.13) and the previous inequality establishes the validity of the left inequality (2.16). Define the operator $\tilde{P} : \Re^{n-1} \to \Re^{n-1}$ defined by

$$\Re^{n-1} \ni (s_2,...,s_n) \to y = (y_2,...,y_n) = \tilde{P}(s) \text{ with } y_i = \max(\lambda, s_i) \quad (6.41)$$

and notice that for every $(s,v) \in \Omega$ it holds that $(\tilde{P}(s), v^* \mathbf{1}_n) \in S$, where $\mathbf{1}_n = (1,1,...,1) \in \mathfrak{R}^n$. Using (2.10), we obtain for all $(s,v) \in \Omega$:

$$dist((s,v), S) \leq \left( \left| \tilde{P}(s) - s \right|^2 + \left| v - v^* \mathbf{1}_n \right|^2 \right)^{1/2} \leq \left( 2H(s,v) + \sum_{i=2}^{n} (\max(\lambda - s_i, 0))^2 \right)^{1/2} \quad (6.42)$$

Since $V$ restricted on the interval $(L, \lambda]$ is a decreasing, continuous function (recall (2.3)), the inverse function of $V$ denoted by $V^{-1} : \mathfrak{R}_+ \to (L, \lambda]$ is a decreasing continuous function with $V^{-1}(0) = \lambda$ and $\lim_{l \to +\infty} (V^{-1}(l)) = L$. It follows that the function $\bar{a}(l) := \lambda - V^{-1}(l)$ defined on $\mathfrak{R}_+$ is a increasing continuous function with $\bar{a}(0) = 0$ and $\lim_{l \to +\infty} (\bar{a}(l)) = \lambda - L$. Moreover, we have $\max(\lambda - s_i, 0) \leq \lambda - V^{-1}(V(s_i)) = \bar{a}(V(s_i)) \leq \bar{a}(H(s,v))$ for all $(s,v) \in \Omega$ and $i = 2,...,n$. Therefore, we obtain from (6.42) for all $(s,v) \in \Omega$:

$$dist((s,v), S) \leq \left( 2H(s,v) + (n-1)(\bar{a}(H(s,v)))^2 \right)^{1/2} \quad (6.43)$$

The existence of a function $a \in K_\infty$ so that $a(dist((s,v), S)) \leq H(s,v)$ for all $(s,v) \in \Omega$ is directly implied by (6.43) (with $a = b^{-1}$ and $b(l) = \left( 2l + (n-1)(\bar{a}(l))^2 \right)^{1/2}$ for $l \geq 0$).

The proof is complete. ◁

**Proof of Theorem 3.2:** Uniqueness of solutions for the initial-value problem (3.8), (3.9) is straightforward and is left to the reader. We prove estimates (3.10), (3.11), and (3.13).

For each $t \geq 0$ we define the function:

$$P_t(r) = \omega r + (1 - \exp(-\omega t)) v_0(r), \text{ for all } r \in \mathfrak{R} \quad (6.44)$$

Due to the facts that $v_0 \in L^\infty(\mathfrak{R})$, $\inf_{x \in \mathfrak{R}} (v'_0(x)) > -\omega$ and $\omega > 0$, it follows that for each $t \geq 0$ the function $P_t : \mathfrak{R} \to \mathfrak{R}$ is increasing with $P_t(\mathfrak{R}) = \mathfrak{R}$. Therefore, for each $t \geq 0$ the inverse function $P_t^{-1} : \mathfrak{R} \to \mathfrak{R}$ is well-defined (of class $C^1(\mathfrak{R})$ and increasing). Moreover, by virtue of the implicit function theorem the function $\mathfrak{R}_+ \times \mathfrak{R} \ni (t,r) \to P_t^{-1}(r)$ is $C^1$ and satisfies for all $(t,r) \in \mathfrak{R}_+ \times \mathfrak{R}$

$$\frac{\partial P_t^{-1}}{\partial r}(r) = \frac{1}{\omega + (1 - \exp(-\omega t)) v'_0(P_t^{-1}(r))}$$

$$\frac{\partial P_t^{-1}}{\partial t}(r) = -\frac{\omega \exp(-\omega t) v_0(P_t^{-1}(r))}{\omega + (1 - \exp(-\omega t)) v'_0(P_t^{-1}(r))} \quad (6.45)$$

It is straightforward to verify that the solution of the initial-value problem (3.8), (3.9) is given by the following formulas:

$$\rho(t,x) = \frac{\omega \rho_0(\xi(t,x))}{\omega + (1 - \exp(-\omega t)) v'_0(\xi(t,x))} \quad (6.46)$$

$$v(t, x) = v^* + \exp(-\omega t)\left(v_0(\xi(t, x)) - v^*\right) \tag{6.47}$$

where

$$\xi(t, x) = P_t^{-1}\left(\omega(x - v^* t) + v^*(1 - \exp(-\omega t))\right) \tag{6.48}$$

Estimates (3.10), (3.11), and (3.12) are direct consequences of the above formulas and the facts that $\inf_{x \in \Re}(v_0'(x)) > -\omega$, $\rho_0(x) > 0$ for all $x \in \Re$.

We next define:

$$\begin{aligned}
\gamma &:= \omega + \inf_{x \in \Re}\left(v_0'(x)\right) \\
K &:= \sup_{x \in \Re}(\rho_0(x)) + \sup_{x \in \Re}(|v_0(x)|) \\
c &:= \sup_{x \in \Re}(|\rho_0'(x)|) + \sup_{x \in \Re}(|v_0''(x)|) \\
L &:= \sup_{x \in \Re}(|v_0'(x)|)
\end{aligned} \tag{6.49}$$

and we notice that equation (6.45) implies the following estimate for all $(t, r) \in \Re_+ \times \Re$:

$$0 < \frac{\partial P_t^{-1}}{\partial r}(r) \leq \frac{1}{\min(\gamma, \omega)} \tag{6.50}$$

Define the function

$$P_\infty(r) = \omega r + v_0(r) \text{, for all } r \in \Re \tag{6.51}$$

and notice that due to the facts that $v_0 \in L^\infty(\Re)$, $\inf_{x \in \Re}(v_0'(x)) > -\omega$ and $\omega > 0$, it follows that for each $t \geq 0$ the function $P_\infty : \Re \to \Re$ is increasing with $P_\infty(\Re) = \Re$. Therefore, for each $t \geq 0$ the inverse function $P_\infty^{-1} : \Re \to \Re$ is well-defined (of class $C^1(\Re)$ and increasing). Moreover, definitions (6.44), (6.51) imply the equation:

$$P_t(r) = P_\infty(r) - \exp(-\omega t)v_0(r) \text{, for all } t \geq 0, r \in \Re \tag{6.52}$$

The above equation in conjunction with (6.49), (6.50) implies the following estimate for all $t \geq 0$, $r \in \Re$:

$$\begin{aligned}
\left|P_t^{-1}(r) - P_\infty^{-1}(r)\right| &= \left|P_t^{-1}(r) - P_t^{-1}\left(P_t\left(P_\infty^{-1}(r)\right)\right)\right| \\
&= \left|P_t^{-1}(r) - P_t^{-1}\left(r - \exp(-\omega t)v_0\left(P_\infty^{-1}(r)\right)\right)\right| \leq \frac{K\exp(-\omega t)}{\min(\gamma, \omega)}
\end{aligned} \tag{6.53}$$

We will show next that estimate (3.13) holds with

$$f(x) := \frac{\omega \rho_0(\zeta(x))}{\omega + v_0'(\zeta(x))} \text{, for all } x \in \Re \tag{6.54}$$

where

$$\zeta(x) = P_\infty^{-1}(\omega x + v^*) \text{, for all } x \in \Re. \tag{6.55}$$

Using the triangle inequality, formulas (6.46), (6.47), definitions (6.49), (6.54) and the facts that $\inf_{x\in\Re}(v_0'(x)) > -\omega$, $\rho_0(x) > 0$ for all $x \in \Re$, we obtain for all $t \geq 0$, $x \in \Re$:

$$\begin{aligned}
&|\rho(t,x) - f(x-v^*t)| \\
&\leq \omega \frac{|\rho_0(\xi(t,x)) - \rho_0(\zeta(x-v^*t))|}{\omega + v_0'(\zeta(x-v^*t))} \\
&\quad + \omega\rho_0(\xi(t,x)) \frac{|v_0'(\zeta(x-v^*t)) - (1-\exp(-\omega t))v_0'(\xi(t,x))|}{(\omega + (1-\exp(-\omega t))v_0'(\xi(t,x)))(\omega + v_0'(\zeta(x-v^*t)))} \\
&\leq \gamma^{-1}\omega|\rho_0(\xi(t,x)) - \rho_0(\zeta(x-v^*t))| \\
&\quad + \frac{\omega K}{\gamma\min(\gamma,\omega)}|v_0'(\zeta(x-v^*t)) - (1-\exp(-\omega t))v_0'(\xi(t,x))| \\
&\leq \gamma^{-1}\omega|\rho_0(\xi(t,x)) - \rho_0(\zeta(x-v^*t))| + \frac{\omega K}{\gamma\min(\gamma,\omega)}|v_0'(\zeta(x-v^*t)) - v_0'(\xi(t,x))| \\
&\quad + \frac{\omega K \exp(-\omega t)}{\gamma\min(\gamma,\omega)}|v_0'(\xi(t,x))| \\
&\leq \frac{\omega c}{\gamma}\left(1 + \frac{K}{\min(\gamma,\omega)}\right)|\xi(t,x) - \zeta(x-v^*t)| + L\frac{\omega K \exp(-\omega t)}{\gamma\min(\gamma,\omega)}
\end{aligned} \quad (6.56)$$

Moreover, using the triangle inequality, (6.48) and definitions (6.49), we get for all $t \geq 0$, $x \in \Re$:

$$\begin{aligned}
|\xi(t,x) - \zeta(x-v^*t)| &= |P_t^{-1}(\omega(x-v^*t) + v^*(1-\exp(-\omega t))) - P_\infty^{-1}(\omega(x-v^*t) + v^*)| \\
&\leq |P_t^{-1}(\omega(x-v^*t) + v^*(1-\exp(-\omega t))) - P_t^{-1}(\omega(x-v^*t) + v^*)| \\
&\quad + |P_t^{-1}(\omega(x-v^*t) + v^*) - P_\infty^{-1}(\omega(x-v^*t) + v^*)| \\
&\leq \frac{v^*\exp(-\omega t)}{\min(\gamma,\omega)} + |P_t^{-1}(\omega(x-v^*t) + v^*) - P_\infty^{-1}(\omega(x-v^*t) + v^*)| \\
&\leq \frac{(K+v^*)\exp(-\omega t)}{\min(\gamma,\omega)}
\end{aligned} \quad (6.57)$$

Combining (6.56) and (6.57) we obtain for all $t \geq 0$, $x \in \Re$:

$$|\rho(t,x) - f(x-v^*t)|\exp(\omega t) \leq \frac{\omega}{\gamma\min(\gamma,\omega)}\left(c(K+v^*)\left(1 + \frac{K}{\min(\gamma,\omega)}\right) + LK\right)$$

Estimate (3.13) is a direct consequence of the above inequality. The proof is complete.  ◁

**Formal Derivation of the Macroscopic Model (3.4)-(3.5):** Assume that each vehicle has mass $\frac{m}{n}$ and length $\frac{\sigma}{n}$. Moreover, define $\lambda = \frac{m}{n\bar{\rho}}$, $L = \frac{m}{n\rho_{\max}}$ for $n = 2, 3, \ldots$ and notice that properties (3.1) guarantee properties (2.3) for the function

$$V(s) := \Phi\left(\frac{m}{ns}\right). \tag{6.58}$$

Consider a solution $(x(t), v(t)) \in \Re^n \times \Re^n$ of the microscopic model (2.1) with $V$ defined by (6.58), $x(t) = (x_1(t), \ldots, x_n(t)) \in \Re^n$ and $v(t) = (v_1(t), \ldots, v_n(t)) \in \Re^n$. Moreover, consider $C^1$ density and speed functions that satisfy the equations

$$\rho(t, x_i(t)) = \frac{m}{n(x_{i-1}(t) - x_i(t))}, \text{ for } t \geq 0, \ i = 2, \ldots, n \tag{6.59}$$

$$v(t, x_i(t)) = v_i(t), \text{ for } t \geq 0, \ i = 1, \ldots, n. \tag{6.60}$$

Using the definition (6.59) above and (2.1), we obtain

$$\frac{d}{dt}\rho(t, x_i(t)) = -\rho(t, x_i(t)) \frac{v_{i-1}(t) - v_i(t)}{x_{i-1}(t) - x_i(t)}, \text{ for } t \geq 0, \ i = 2, \ldots, n. \tag{6.61}$$

Using the chain rule, we also have from (6.59), (6.60) and (2.1) that

$$\frac{d}{dt}\rho(t, x_i(t)) = \rho_t(t, x_i(t)) + \rho_x(t, x_i(t))v(t, x_i(t)), \text{ for } t \geq 0, \ i = 2, \ldots, n. \tag{6.62}$$

Thus, we get from (6.60), (6.61) and (6.62) for all $t \geq 0$, $i = 2, \ldots, n$:

$$\rho_t(t, x_i(t)) + \rho_x(t, x_i(t))v(t, x_i(t)) + \rho(t, x_i(t))\frac{v(t, x_{i-1}(t)) - v(t, x_i(t))}{x_{i-1}(t) - x_i(t)} = 0. \tag{6.63}$$

Assuming that $x_{i-1}(t) - x_i(t) \to 0$ as $n \to +\infty$, we have $\frac{v(t, x_{i-1}(t)) - v(t, x_i(t))}{x_{i-1}(t) - x_i(t)} \to v_x(t, x_i(t))$ for $i = 2, \ldots, n$. Consequently, we obtain from (6.63) the continuity equation

$$\rho_t(t, x) + \rho_x(t, x)v(t, x) + \rho(t, x)v_x(t, x) = 0, \text{ for } t \geq 0, \ x \in \Re. \tag{6.64}$$

Next, using the definition of speed (6.60) and (2.1), (2.3), we obtain

$$\frac{d}{dt}v(t, x_i(t)) = -\left(\mu + g\left(V'(s_i(t)) - V'(s_{i+1}(t))\right)\right)(v(t, x_i(t)) - v^*) + V'(s_i(t)) - V'(s_{i+1}(t)),$$
$$\text{for } t \geq 0, \ i = 2, \ldots, n-1.$$

We also obtain from the chain rule

$$\frac{d}{dt}v(t,x_i(t)) = v_t(t,x_i(t)) + v_x(t,x_i(t))v(t,x_i(t)), \text{ for } t \geq 0, \ i = 2,\ldots,n-1.$$

The equations above imply that

$$v_t(t,x_i(t)) + v(t,x_i(t))v_x(t,x_i(t)) = -\big(\mu + g(\Xi_i(t))\big)(v(t,x_i(t)) - v^*) + \Xi_i(t),$$
$$\text{for } t \geq 0, \ i = 2,\ldots,n-1 \qquad (6.65)$$

where

$$\Xi_i(t) := V'(s_i(t)) - V'(s_{i+1}(t)), \text{ for } t \geq 0, \ i = 2,\ldots,n-1. \qquad (6.66)$$

From definitions (6.58), (6.59), (6.66) and Taylor's Theorem we have that

$$\begin{aligned}
\Xi_i(t) &= V'(s_i(t)) - V'(s_{i+1}(t)) \\
&= V'\!\left(\frac{m}{n\rho(t,x_i(t))}\right) - V'\!\left(\frac{m}{n\rho(t,x_{i+1}(t))}\right) \\
&= -\frac{n}{m}\Big((\rho(t,x_i(t)))^2 \Phi'(\rho(t,x_i(t))) - (\rho(t,x_{i+1}(t)))^2 \Phi'(\rho(t,x_{i+1}(t)))\Big) \\
&\approx -\frac{n}{m}\big(\rho(t,x_i(t)) - \rho(t,x_{i+1}(t))\big) \\
&\quad \times \Big(2\rho(t,x_i(t))\Phi'(\rho(t,x_i(t))) + (\rho(t,x_{i+1}(t)))^2 \Phi''(\rho(t,x_i(t)))\Big) \\
&= -\frac{n}{m}\big(x_i(t) - x_{i+1}(t)\big)\frac{\rho(t,x_i(t)) - \rho(t,x_{i+1}(t))}{x_i - x_{i+1}} \\
&\quad \times \Big(2\rho(t,x_i(t))\Phi'(\rho(t,x_i(t))) + (\rho(t,x_i(t)))^2 \Phi''(\rho(t,x_i(t)))\Big)
\end{aligned} \qquad (6.67)$$

Assuming that $x_{i-1}(t) - x_i(t) \to 0$ as $n \to +\infty$, we have $\dfrac{\rho(t,x_i(t)) - \rho(t,x_{i+1}(t))}{x_i(t) - x_{i+1}(t)} \to \rho_x(t,x_i(t))$ for $i = 2,\ldots,n$. Consequently, by defining the function $\Xi(t,x)$ that satisfies $\Xi(t,x_i(t)) = \Xi_i(t)$ for all $t \geq 0$ and $i = 2,\ldots,n-1$, we conclude from (6.65), (6.67), and definition (6.59) that equations (3.4), (3.5) hold.

## 7. Conclusions

The paper introduces a new bidirectional microscopic inviscid Adaptive Cruise Control (ACC) model that uses only spacing information from the preceding and following vehicles. *KL* estimates that guarantee uniform convergence properties of the ACC model to the set of equilibria are provided. Moreover, the corresponding macroscopic model is derived, consisting of the continuity equation and a momentum equation that contains a highly nonlinear relaxation term. It is shown that, if the density is sufficiently small, then the solution of the macroscopic model approaches the equilibrium speed (in the sup norm); while the density converges exponentially to a traveling wave. Numerical simulations are also provided to illustrate the properties of the microscopic and macroscopic inviscid models.

Future work will involve the study of the movement of vehicles under different cruise controllers. The cruise controller that was applied in the present work does not give rise to a viscous term (and that is why is termed "inviscid") and is based on the classical mechanics of particles (vehicles) with a particularly strong "friction-like" term (the relaxation term) that guarantees that the speed remains within specific bounds (between 0 and $v_{max}$). Other cruise controllers are envisaged as follows:

- Cruise controllers may be based on "pseudo-relativistic mechanics" of particles (vehicles). In relativistic mechanics, the speeds of the particles are, in absolute value, always less than the speed of light and this is guaranteed by accelerations of the form $\dot{v} \approx -g(|v|^2)V'(x)$ (when the motion is on a straight line), where $g$ is a specific function that tends to zero when the absolute value of the speed of the body tends to the speed of light. Similarly, for vehicles one can consider cruise controllers that are of the form $\dot{v}_i = g(v_i)(V'(s_i) - V'(s_{i+1}) - \mu(v_i - v^*))$, where $\mu > 0$ is the relaxation constant and $g$ is a specific function that tends to zero when the speed of the body tends to 0 or $v_{max}$. The produced model can be termed as a "bidirectional, microscopic, pseudo-relativistic ACC model" and will feature important differences from the bidirectional, microscopic, inviscid ACC model, since the controller gains of the pseudo-relativistic ACC will be much lower than the gains of the non-relativistic (classical) model that was studied in the present paper. The much lower gains will allow for higher sampling periods in the case of the sampled-data implementation of the ACC, but will also induce much slower convergence to the set of equilibrium points.

- Cruise controllers may include a viscous term, i.e. a term of the form $\kappa(s_i)(v_{i-1} - v_i) + \kappa(s_{i+1})(v_{i+1} - v_i)$ in the acceleration equation, where $\kappa$ is a specific non-negative function. This term considers the relative speed of subsequent vehicles and gives rise to a viscous term in the macroscopic model; this can be considered in either a "pseudo-relativistic" setting or a "classical (Newtonian)" setting. The produced model can be termed as a "bidirectional microscopic viscous ACC model" and will again feature important differences from the bidirectional, microscopic, inviscid ACC model, since the viscous term enhances energy dissipation.

The corresponding macroscopic models of the above cruise controllers will also be studied. The macroscopic models will give rise to systems of PDEs, which are similar to the equations describing the flow of a compressible (viscous or inviscid, relativistic or not) fluid. It is interesting to highlight that by selecting an appropriate cruise controller, we can design a corresponding "artificial fluid" [30], in the sense that we have the possibility to impose certain desired properties to the traffic fluid. This is important in the era of connected and automated vehicles, because, while safety and passenger convenience issues are taken care of at the microscopic level, it is the macroscopic level that reflects important emerging traffic flow features, such as flow level, capacity and stability. Thus, the present work gives only a glimpse of the future of cruise controller design.

## Acknowledgments

The research leading to these results has received funding from the European Research Council under the European Union's Horizon 2020 Research and Innovation programme/ ERC Grant Agreement n. [833915], project TrafficFluid.